\documentclass{article}

\usepackage{multicol}
\usepackage{graphicx}
\usepackage{amsfonts}
\usepackage{amsmath}
\usepackage{url}
\usepackage{color}
\usepackage[margin=1.0in]{geometry}

\title{Stability-Optimized High Order Methods and Stiffness Detection for
	Pathwise Stiff Stochastic Differential Equations\thanks{
			This work was partially supported by NIH grants R01GM107264 and P50GM76516
			and NSF grants DMS1562176 and DMS1161621. This material is based upon
			work supported by the National Science Foundation Graduate Research
			Fellowship under Grant No. DGE-1321846, the National Academies of
			Science, Engineering, and Medicine via the Ford Foundation, and the
			National Institutes of Health Award T32 EB009418. Its contents are
			solely the responsibility of the authors and do not necessarily represent
			the official views of the NIH.}}

\author{Christopher Rackauckas\thanks{Department of Mathematics, University of California, Irvine, Irvine,
		CA 92697, USA and Center for Complex Biological Systems, University
		of California, Irvine, Irvine, CA 92697, USA (contact@chrisrackauckas.com).}
	\and Qing Nie\thanks{Department of Mathematics, University of California, Irvine, Irvine,
		CA 92697, USA and Center for Complex Biological Systems, University
		of California, Irvine, Irvine, CA 92697, USA and Department of Developmental
		and Cell Biology, University of California, Irvine, Irvine, CA 92697,
		USA (qnie@math.uci.edu).}}

\begin{document}
	
\maketitle

\begin{abstract}
	Stochastic differential equations (SDE) often exhibit large random
	transitions. This property, which we denote as pathwise stiffness,
	causes transient bursts of stiffness which limit the allowed step
	size for common fixed time step explicit and drift-implicit integrators.
	We present four separate methods to efficiently handle this stiffness.
	First, we utilize a computational technique to derive stability-optimized
	adaptive methods of strong order 1.5 for SDEs. The resulting explicit
	methods are shown to exhibit substantially enlarged stability regions
	which allows for them to solve pathwise stiff biological models orders
	of magnitude more efficiently than previous methods like SRIW1 and
	Euler-Maruyama. Secondly, these integrators include a stiffness estimator
	which allows for automatically switching between implicit and explicit
	schemes based on the current stiffness. In addition, adaptive L-stable
	strong order 1.5 implicit integrators for SDEs and stochastic differential
	algebraic equations (SDAEs) in mass-matrix form with additive noise
	are derived and are demonstrated as more efficient than the explicit
	methods on stiff chemical reaction networks by nearly 8x. Lastly,
	we developed an adaptive implicit-explicit (IMEX) integration method
	based off of a common method for diffusion-reaction-convection PDEs
	and show numerically that it can achieve strong order 1.5. These methods
	are benchmarked on a range of problems varying from non-stiff to extreme
	pathwise stiff and demonstrate speedups between 5x-6000x while showing
	computationally infeasibility of fixed time step integrators on many
	of these test equations.
\end{abstract}

\section{Introduction}

Stochastic differential equations (SDEs) are dynamic equations of
the form
\begin{equation}
dX_{t}=f(t,X_{t})dt+g(t,X_{t})dW_{t},
\end{equation}
where $X_{t}$ is a $d$-dimensional vector, $f:\mathbb{R}^{d}\rightarrow\mathbb{R}^{d}$
is the drift coefficient, and $g:\mathbb{R}^{d}\rightarrow\mathbb{R}^{d\times m}$
is the diffusion coefficient which describes the amount and mixtures
of the noise process $W_{t}$ which is a $m$-dimensional Brownian
motion. SDEs are of interest in scientific disciplines because they
can exhibit behaviors which are not found in deterministic models.
For example, An ODE model of a chemical reaction network may stay
at a constant steady state, but in the presence of randomness the
trajectories may be switching between various steady states \cite{RN3805,RN3814,RN3358}.
In many cases, these unique features of stochastic models are pathwise-dependent
and are thus not a property of the evolution of the mean trajectory.
However, these same effects cause random events of high numerical
stiffness, which we denote as pathwise stiffness, which can cause
difficulties for numerical integration methods.

\begin{figure}
	\begin{centering}
		\includegraphics[scale=0.5]{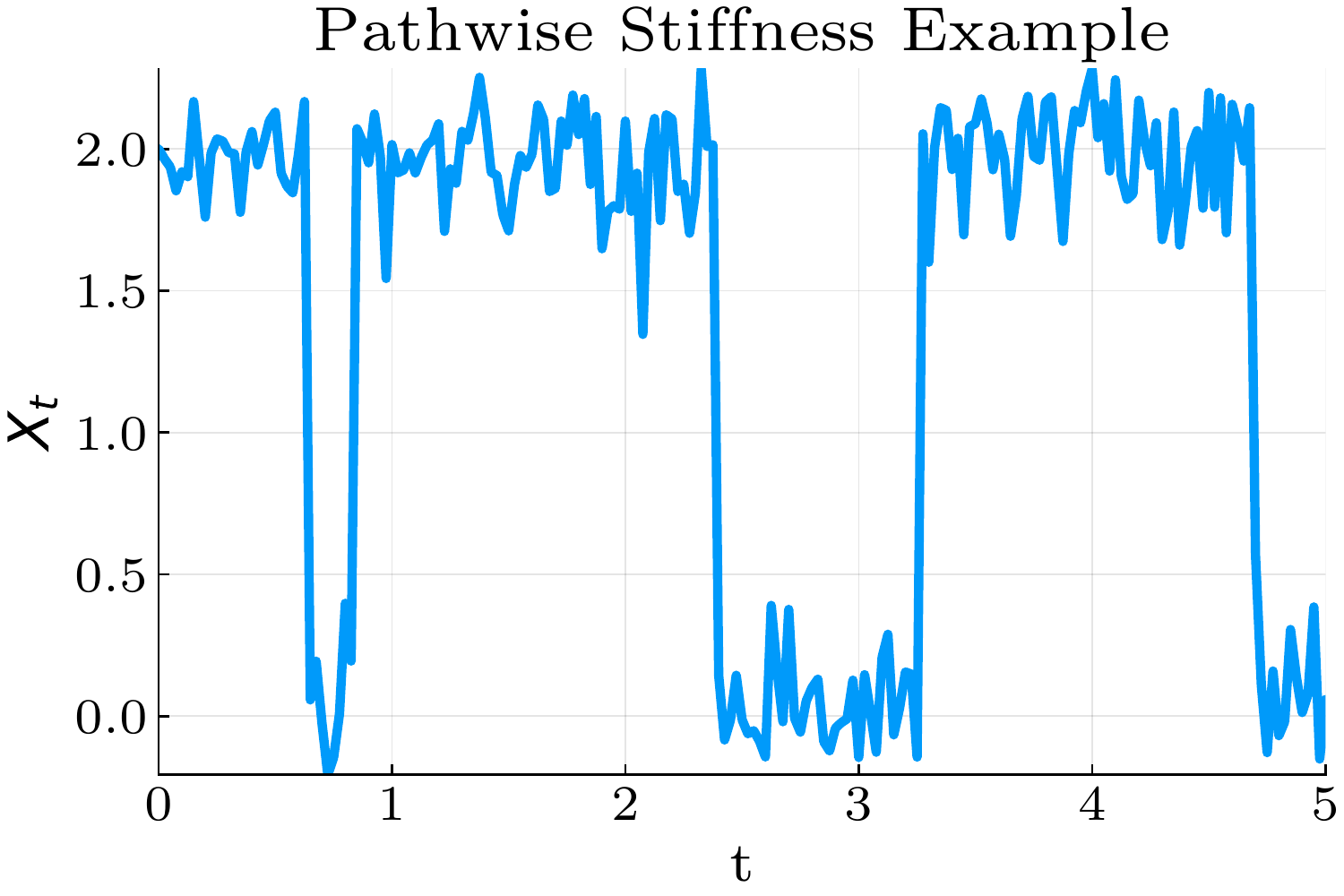}
		\par\end{centering}
	\caption{\textbf{Example of a Pathwise Stiff Solution}. Depicted is a sample
		trajectory of Equation \ref{eq:pathwise_stiff} solved using the SOSRI
		methods developed in this manuscript with $reltol=abstol=10^{-2}$.
		\label{fig:pathwise_stiff}}
\end{figure}

A minimal example of pathwise stiffness is demonstrated in the equation
\begin{equation}
dX_{t}=\left[-1000X_{t}\left(1-X_{t}\right)\left(2-X_{t}\right)\right]dt+g(t,X_{t})dW_{t},\thinspace\thinspace\thinspace X_{0}=2,\thinspace\thinspace\thinspace t\in\left[0,5\right].\label{eq:pathwise_stiff}
\end{equation}
with additive noise $g(t,X_{t})=10$ where a sample trajectory is
shown in Figure \ref{fig:pathwise_stiff}. This equation has two stable
steady states, one at $X=0$ and another at $X=2$, which the solution
switches between when the noise is sufficiently large. While near
a steady state the derivative is approximately zero making the problem
non-stiff, during these transitions the derivative of the drift term
reaches a maximum of $\approx400$. This means that in order to be
stable, explicit Stochastic Runge-Kutta (SRK) must have a small $\Delta t$.
This display of large, transient, and random switching behavior in
a given trajectory causes stochastic bursts of numerical stiffness,
a phenomena which we will denote pathwise stiffness. The fixed time
step Euler-Maruyama method would require $dt<4\times10^{-3}$ to be
stable for most trajectories, thus requiring greater than $2\times10^{4}$
steps to solve this 1-dimensional SDE. In many cases the switching
behavior can be rare (due to smaller amounts of noise) or can happen
finitely many times like in the multiplicative noise version with
$g(t,X_{t})=10X_{t}$. Yet even if these switches are only a small
portion of the total time, the stability requirement imposed by their
existence determines the possible stepsizes and thus has a large contribution
to the overall computational cost. While implicit methods can be used
to increase the stability range, this can vastly increase the overall
computational cost of each step, especially in the case large systems
of SDEs like discretizations of stochastic reaction-diffusion equations.
In addition, implicit solvers have in practice a smaller stability
region due to requiring convergence of the quasi-Newton solvers for
the implicit steps. This problem is mitigated in ODE software by high-quality
stage predictors given by extrapolation algorithms for good initial
conditions for the Newton steps \cite{RN3790}. However, there are
no known algorithms for stage predictors in the presence of large
noise bursts and thus we will demonstrate that classic implicit solvers
have a form of instability. Thus both fixed time step explicit and
implicit solvers are inadequate for efficiently handling this common
class of SDEs.

Since these features exist in the single trajectories of the random
processes, methods which attempt to account for the presence of such
bursts must do so on each individual trajectory in order to be efficient.
In previous work, the authors have shown that by using adaptive time-stepping,
a stochastic reaction network of 19 reactants is able to be solved
with an average time step 100,000 times larger than the value that
was found necessary for stability during the random stiff events for
a high order SRK method \cite{RN3787}. This demonstrated that the
key to solving these equations efficiently required controlling the
time steps in a pathwise manner. However, the methods were still largely
stability-bound, meaning the chosen tolerances to solve the model
were determined by what was necessary for stability and was far below
the error necessary for the application. The purpose of this investigation
is to develop numerical methods with the ability to better handle
pathwise stiffness and allow for efficient solving of large Monte
Carlo experiments.

We approach this problem through four means. First, we develop adaptive
stability-optimized SRK methods with enlarged stability regions. This
builds off of similar work for ODE integrators which optimize the
coefficients of a Butcher tableau to give enhanced stability \cite{RN3514,RN3515,RN3517}.
Similar to the Runge-Kutta Chebyschev methods \cite{RN3790} (and
the S-ROCK extension to the stochastic case \cite{RN3799,RN3800,RN3801}),
these methods are designed to be efficient for equations which display
stiffness without fully committing to implicit solvers. Given the
complexity of the stochastic stability equations and order conditions,
we develop a novel and scalable mechanism for the derivation of ``optimal''
Runge-Kutta methods. We use this method to design stability-optimized
methods for additive noise and diagonal noise SDEs. We show through
computational experiments that these adaptive stability-optimized
SRK methods can adequately solve transiently stiff equations without
losing efficiency in non-stiff problems.

On the other hand, to handle extreme stiffness we develop implicit
RK methods for SDEs and stochastic differential algebraic equations
(SDAEs) in mass matrix form with additive noise. We extend the definition
of L-stability to additive noise SDEs and develop two strong order
1.5 methods: a fully implicit 2-stage L-stable method and an extension
of the a well-known L-stable explicit first stage singly diagonally
implicit RK (ESDIRK) method due to Kennedy and Carpenter which is
commonly used for convection-diffusion-reaction equations \cite{RN3792}.
To the author's knowledge, these are the first high order adaptive
L-stable methods for SDEs and the first adaptive proposed SDAE integrators.
In addition, to extend the utility of these additive noise methods,
we derive an extension of the methods for additive SDEs to affine
SDEs (mixed multiplicative and additive noise terms) through a Lamperti
transformation \cite{RN3525}. Lastly, in order to handle extreme
transient stiffness, for each of these types of methods we derive
computationally cheap methods for detecting stiffness and switching
between implicit and explicit integrators in the presence of stiffness.
We show that these methods can robustly detect pathwise stiff transients
and thus can serve as the basis for automatic switching methods for
SDEs. Together we test on non-stiff, semi-stiff, and stiff equations
with 2 to $6\times20\times100$ SDEs from biological literature and
show speedups between 6x-60x over the previous adaptive SRIW1 algorithm,
and demonstrate the infeasibility of common explicit and implicit
methods (Euler-Maruyama, Runge-Kutta Milstein, Drift-Implicit Stochastic
$\theta$-Method, and Drift-Implicit $\theta$ Runge-Kutta Milstein)
found as the basis of many SDE solver packages \cite{RN3803,RN3802,RN3804}.

\section{Adaptive Strong Order 1.0/1.5 SRK Methods for Additive and Diagonal
	Noise SDEs}

The class of methods we wish to study are the adaptive strong order
1.5 SRK methods for diagonal noise \cite{RN2707,RN3787}. Diagonal
noise is the case where the diffusion term $g$ is diagonal matrix
$\left(\sigma_{i}X_{t}^{i}\right)$ and includes phenomenological
noise models like multiplicative and affine noise. The diagonal noise
methods utilize the same general form and order conditions as the
methods for scalar noise so we use their notation for simplicity.
The strong order 1.5 methods for scalar noise are of the form

\begin{align}
X_{n+1} & =X_{n}+\sum_{i=1}^{s}\alpha_{i}f\left(t_{n}+c_{i}^{(0)}h,H_{i}^{(0)}\right)+\\
        & \sum_{i=1}^{s}\left(\beta_{i}^{(1)}I_{(1)}+\beta_{i}^{(2)}\frac{I_{(1,1)}}{\sqrt{h}}+\beta_{i}^{(3)}\frac{I_{(1,0)}}{h}+\beta_{i}^{(4)}\frac{I_{(1,1,1)}}{h}\right)g\left(t_{n}+c_{i}^{(1)}h\right)\label{eq:update}
\end{align}
with stages
\begin{align}
H_{i}^{(0)} & =X_{n}+\sum_{j=1}^{s}A_{ij}^{(0)}f\left(t_{n}+c_{j}^{(0)}h,H_{j}^{(0)}\right)h+\sum_{j=1}^{s}B_{ij}^{(0)}g\left(t_{n}+c_{j}^{(1)}h,H_{j}^{(1)}\right)\frac{I_{(1,0)}}{h}\label{eq:stages}\\
H_{i}^{(1)} & =X_{n}+\sum_{j=1}^{s}A_{ij}^{(1)}f\left(t_{n}+c_{j}^{(0)}h,H_{j}^{(0)}\right)h+\sum_{j=1}^{s}B_{ij}^{(1)}g\left(t_{n}+c_{j}^{(1)}h,H_{j}^{(1)}\right)\sqrt{h}\nonumber
\end{align}
where the $I_{j}$ are the Wiktorsson approximations to the iterated
stochastic integrals \cite{RN3175}. In the case of additive noise,
defined as having the diffusion coefficient satisfy $g(t,X_{t})\equiv g(t)$
, reduces to the form

\begin{equation}
X_{n+1}=X_{n}+\sum_{i=1}^{s}\alpha_{i}f\left(t_{n}+c_{i}^{(0)}h,H_{i}^{(0)}\right)+\sum_{i=1}^{s}\left(\beta_{i}^{(1)}I_{(1)}+\beta_{i}^{(2)}\frac{I_{(1,0)}}{h}\right)g\left(t_{n}+c_{i}^{(1)}h\right)\label{eq:add_update}
\end{equation}
with stages
\begin{equation}
H_{i}^{(0)}=X_{n}+\sum_{j=1}^{s}A_{ij}^{(0)}f\left(t_{n}+c_{j}^{(0)}h,H_{j}^{(0)}\right)h+\sum_{j=1}^{s}B_{ij}^{(0)}g\left(t_{n}+c_{j}^{(1)}h\right)\frac{I_{(1,0)}}{h}.\label{eq:add_stages}
\end{equation}
The tuple of coefficients $\left(A^{(j)},B^{(j)},\beta^{(j)},\alpha\right)$
thus fully determines the SRK method. These coefficients must satisfy
the constraint equations described in Appendix \ref{subsec:Order-Conditions-for-SRI}
in order to receive strong order 1.5. These methods are appended with
error estimates

\begin{eqnarray*}
	E_{D} & = & \left|\Delta t\sum_{i\in I_{1}}(-1)^{\sigma(i)}f\left(t_{n}+c_{i}^{(0)}\Delta t,H_{i}^{(0)}\right)\right|\mbox{or }E_{D}=\Delta t\sum_{i\in I_{1}}\left|f\left(t_{n}+c_{i}^{(0)}\Delta t,H_{i}^{(0)}\right)\right|\\
	E_{N} & = & \left|\sum_{i\in I_{2}}\left(\beta_{i}^{(3)}\frac{I_{(1,0)}}{\Delta t}+\beta_{i}^{(4)}\frac{I_{(1,1,1)}}{\Delta t}\right)g\left(t_{n}+c_{i}^{(1)}\Delta t,H_{i}^{(1)}\right)\right|
\end{eqnarray*}
and the rejection sampling with memory (RSwM) algorithm to give it
fully adaptive time-stepping \cite{RN3787}. Thus unlike in the theory
of ordinary differential equations \cite{RN3512,RN3513,RN3519,RN3523,RN3518},
the choice of coefficients for SRK methods does not require explicitly
finding an embedded method when developing an adaptive SRK method
and we will therefore take for granted that each of the derived methods
is adaptive.

\section{Optimized-Stability High Order SRK Methods with Additive Noise}

We use a previous definition of a discrete approximation as numerically
stable if for any finite time interval $\left[t_{0},T\right]$, there
exists a positive constant $\Delta_{0}$ such that for each $\epsilon>0$
and each $\delta\in\left(0,\Delta_{0}\right)$
\begin{equation}
\lim_{\left|X_{0}^{\delta}-\bar{X}_{0}^{\delta}\right|\rightarrow0}\sup_{t_{0}\leq t\leq T}P\left(\left|X_{t}^{\delta}-\bar{X}_{t}^{\delta}\right|\geq\epsilon\right)=0\label{eq:add_stab}
\end{equation}
where $X_{n}^{\delta}$ is a discrete time approximation with maximum
step size $\delta>0$ starting at $X_{0}^{\delta}$ and $\bar{X}_{n}^{\delta}$
respectively starting at $\bar{X}_{n}^{\delta}$ \cite{RN3169}. For
additive noise, we consider the complex-valued linear test equations
\begin{equation}
dX_{t}=\mu X_{t}dt+dW_{t}\label{eq:add_test}
\end{equation}
where $\mu$ is a complex number. In this framework, a scheme which
can be written in the form
\begin{equation}
X_{n+1}^{h}=X_{n}^{h}G\left(\mu h\right)+Z_{n}^{\delta}\label{eq:add_test_2}
\end{equation}
with a constant step size $\delta\equiv h$ and $Z_{n}^{\delta}$
are random variables which do not depend on the $Y_{n}^{\delta}$,
then the region of absolute stability is the set where for $z=\mu h$,
$\left|G(z)\right|<1$.

The additive SRK method can be written as
\begin{equation}
X_{n+1}^{h}=X_{n}^{h}+z\left(\alpha\cdot H^{(0)}\right)+\beta^{(1)}\sigma I_{(1)}+\sigma\beta^{(2)}\frac{I_{(1,0)}}{h}\label{eq:SRK_stab_step1}
\end{equation}
where
\begin{equation}
H^{(0)}=\left(I-zA^{(0)}\right)^{-1}\left(\hat{X_{n}^{h}}+B^{(0)}e\sigma\frac{I_{(1,0)}}{h}\right)\label{eq:SRK_stab_step}
\end{equation}
where $\hat{X_{n}^{h}}$ is the size $s$ constant vector of elements
$X_{n}^{h}$ and $e=\left(1,1,1,1\right)^{T}$. By substitution we
receive
\begin{align}
X_{n+1}^{h} & =X_{n}^{h}\left(1+z\left(\alpha\cdot\left(I-zA^{(0)}\right)^{-1}\right)\right)+\\
& \left(I-zA^{(0)}\right)^{-1}B^{(0)}e\sigma\frac{I_{(1,0)}}{h}+\beta^{(1)}\sigma I_{(1)}+\sigma\beta^{(2)}\frac{I_{(1,0)}}{h}\label{eq:SRK_stab_sub}
\end{align}
This set of equations decouples since the iterated stochastic integral
approximation $I_{j}$ are random numbers and are independent of the
$X_{n}^{h}$. Thus the stability condition is determined by the equation
\begin{equation}
G(z)=1+z\alpha\cdot\left(I-zA^{(0)}\right)^{-1}\label{eq:SRA_stab}
\end{equation}
which one may notice is the stability equation of the drift tableau
applied to a deterministic ODE \cite{RN3516}. Thus the stability
properties of the deterministic Runge-Kutta methods carry over to
the additive noise SRA methods on this test equation. However, most
two-stage tableaus from ODE research were developed to satisfy higher
order ODE order constraints which do not apply. Thus we will instead
look to maximize stability while satisfying the stochastic order constraints.

\subsection{Explicit Methods for Non-Stiff SDEs with Additive Noise}

\subsubsection{Stability-Optimal 2-Stage Explicit SRA Methods}

For explicit \\
methods, $A^{(0)}$ and $B^{(0)}$ are lower diagonal
and we receive the simplified stability function
\begin{equation}
G(z)=1+A_{21}z^{2}\alpha_{2}+z\left(\alpha_{1}+\alpha_{2}\right)\label{eq:SRA_stab_2}
\end{equation}
for a two-stage additive noise SRK method. For this method we will
find the method which optimizes the stability in the real part of
$z$. Thus we wish to find $A^{(0)}$ and $\alpha$ s.t. the negative
real roots of $\left|G(z)\right|=1$ are minimized. By the quadratic
equation we see that there exists only a single negative root: $z=\frac{1-\sqrt{1+8\alpha_{2}}}{2\alpha_{2}}$.
Using Mathematica's minimum function, we determine that the minimum
value for this root subject to the order constraints is $z=\frac{3}{4}\left(1-\sqrt{\frac{19}{3}}\right)\approx-1.13746$.
This is achieved when $\alpha=\frac{2}{3}$, meaning that the SRA1
method due to Rossler achieves the maximum stability criteria. However,
given extra degrees of freedom, we attempted to impose that $c_{1}^{(0)}=c_{1}^{(1)}=0$
and $c_{2}^{(0)}=c_{2}^{(1)}=1$ so that the error estimator spans
the whole interval. This can lead to improved robustness of the adaptive
error estimator. In fact, when trying to optimize the error estimator's
span we find that there is no error estimator which satisfies $c_{2}^{(0)}>\frac{3}{4}$
which is the span of the SRA1 method \cite{RN2707}. Thus SRA1 is
the stability-optimized 2-stage explicit method which achieves the
most robust error estimator.

\begin{align}
A^{(0)} & =\left(\begin{array}{cc}
0 & 0\\
\frac{3}{4} & 0
\end{array}\right),\thinspace\thinspace\thinspace B^{(0)}=\left(\begin{array}{cc}
0 & 0\\
\frac{3}{2} & 0
\end{array}\right),\thinspace\thinspace\thinspace\thinspace\alpha=\left(\begin{array}{c}
\frac{1}{3}\\
\frac{2}{3}
\end{array}\right)\nonumber \\
\beta^{(1)} & =\left(\begin{array}{c}
1\\
0
\end{array}\right),\thinspace\thinspace\thinspace\beta^{(2)}=\left(\begin{array}{c}
-1\\
1
\end{array}\right),\thinspace\thinspace\thinspace c^{(0)}=\left(\begin{array}{c}
0\\
\frac{3}{4}
\end{array}\right),\,\,\,c^{(1)}=\left(\begin{array}{c}
1\\
0
\end{array}\right)\label{eq:SRA1}
\end{align}

\subsubsection{Stability-Optimal 3-Stage Explicit SRA Methods}

For the 3-stage SRA method, we receive the simplified stability function
\begin{equation}
G(z)=A_{21}A_{31}\alpha_{3}z^{3}+A_{21}\alpha_{2}z^{2}+A_{31}\alpha_{3}z^{2}+A_{32}\alpha_{3}z^{2}+\alpha_{1}z+\alpha_{2}z+\alpha_{3}z+1\label{eq:SRA_stab_3}
\end{equation}
To optimize this method, we attempted to use the same techniques
as before and optimize the real values of the negative roots. However,
in this case we have a cubic polynomial and the root equations are
more difficult. Instead, we turn to a more general technique to handle
the stability optimization which will be employed in later sections
as well. To do so, we generate an optimization problem which we can
numerically solve for the coefficients. To simplify the problem, we
let $z\in\mathbb{R}$ and define the function:

\begin{equation}
f\left(z,w;N,M\right)=\int_{D}\chi_{G(z)\leq1}(z)dz\label{eq:SRA_integral}
\end{equation}
Notice that $f$ is the area of the stability region when $D$ is
sufficiently large. Thus we define the stability-optimized SRK method
for additive noise SDEs as the set of coefficients which achieves

\begin{align}
\max_{A^{(i)},B^{(i)},\beta^{(i)},\alpha} & f(z)\label{eq:SRA_maximize}\\
\text{subject to: } & \text{Order Constraints}\nonumber
\end{align}
In all cases we impose $0<c_{i}^{(0)},c_{i}^{(1)}<1$. We use the
order constraints to simplify the problem to a nonlinear optimization
problem on 14 variables with 3 equality constraints and 4 inequality
constraints (with bound constraints on the 10 variables). However,
we found that simplifying the problem even more to require $c_{1}^{(0)}=c_{1}^{(1)}=0$
and $c_{3}^{(0)}=c_{3}^{(1)}=1$ did not significantly impact the
stability regions but helps the error estimator and thus we reduced
the problem to 10 variables, 3 equality constraints, and 2 inequality
constraints. This was optimized using the COBYLA local optimization
algorithm \cite{RN3788,RN3789} with randomized initial conditions
100 times and all gave similar results. In the Mathematica notebook
we show the effect of changing the numerical integration region $D$
on the results, but conclude that a $D$ which does not bias the result
for better/worse real/complex handling does not improve the result.
The resulting algorithm, SOSRA, we given by the coefficients in table
in Section \ref{subsec:SOSRA}. Lastly, we used the condition that
$c_{2}^{(0)}=c_{3}^{(0)}=c_{2}^{(1)}=c_{3}^{(1)}=1$ to allow for
free stability detection (discussed in Section \ref{subsec:StabilityDetection}).
The method generated with this extra constraint is SOSRA2 whose coefficients
are in the table in Section \ref{subsec:SOSRA2}. These methods have
their stability regions compared to SRA1 and SRA3 in Figure \ref{fig:SOSRA-Stability-Regions.}
where it is shown that the SOSRA methods more than doubles the allowed
time steps when the eigenvalues of the Jacobian are dominated by the
real part.
\begin{center}
	\begin{figure}
		\begin{centering}
			\includegraphics[scale=0.5]{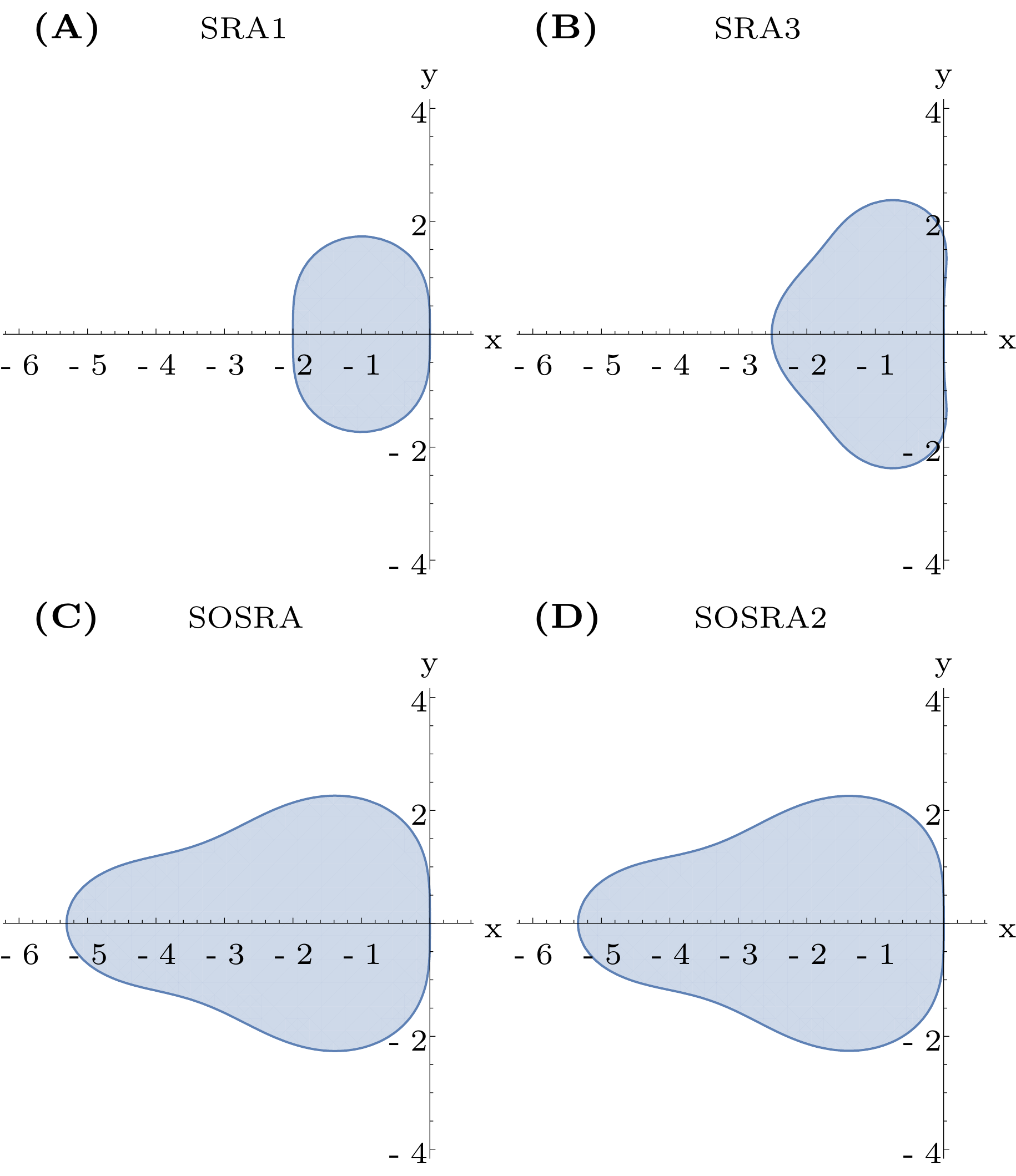}
			\par\end{centering}
		\caption{\textbf{SOSRA Stability Regions}. The stability regions ($\left|G(z)\right|<1$)
			are plotted in the $(x,y)$-plane for $z=x+iy$. \textbf{(A)} SRA1.
			\textbf{(B) }SRA3. \textbf{(C) }SOSRA. \textbf{(D)} SOSRA2\label{fig:SOSRA-Stability-Regions.}}
	\end{figure}
	\par\end{center}

\subsection{Drift Implicit Methods for Stiff SDEs with Additive Noise}

\subsubsection{An L-Stable 2-Stage (Semi-)Implicit SRA Method}

It's clear that, as in the case for deterministic equations, the explicit
methods cannot be made A-stable. However, the implicit two-stage additive
noise SRK method is determined by
\begin{equation}
\small
G(z)=\frac{z(A_{11}(A_{22}z-\alpha_{2}z-1)+A_{12}z(\alpha_{1}-A_{21})+A_{21}\ensuremath{\alpha_{2}}z-A_{22}(\alpha_{1}z+1)+\alpha_{1}+\alpha_{2})+1}{A_{11}z(A_{22}z-1)-z(A_{12}A_{21}z+A_{22})+1}\label{eq:SRA_stab2_implicit}
\end{equation}
which is $A$-stable if
\begin{align}
A_{11}z(A_{22}z-1)-z(A_{12}A_{21}z+A_{22})+1&>z(A_{11}(A_{22}z-\alpha_{2}z-1)+A_{12}z(\alpha_{1}-A_{21})\\
&+A_{21}\ensuremath{\alpha_{2}}z-A_{22}(\alpha_{1}z+1)+\alpha_{1}+\alpha_{2})+1.\label{eq:SRA_stab2_ineq}
\end{align}
Notice that the numerator equals the denominator if and only if $z=0$
or
\begin{equation}
z=\frac{\alpha_{1}+\alpha_{2}}{\left(A_{22}-A_{12}\right)\alpha_{1}+\left(A_{11}-A_{21}\right)\alpha_{2}}.\label{eq:SRA_stab2_z}
\end{equation}
From the order conditions we know that $\alpha_{1}+\alpha_{2}=1$
which means that no root exists with $Re(z)<0$ if $\left(A_{22}-A_{12}\right)\alpha_{1}+\left(A_{11}-A_{21}\right)\alpha_{2}>0$.
Thus under these no roots conditions, we can determine A-stability
by checking the inequality at $z=1$, which gives $1>\left(A_{22}-A_{12}\right)\alpha_{1}+\left(A_{11}-A_{21}\right)\alpha_{2}$.
Using the order condition, we have a total of four constraints on
the $A^{(0)}$ and $\alpha$:

\begin{align}
\left(A_{11}+A_{12}\right)\alpha_{1}+\left(A_{21}+A_{22}\right)\alpha_{2} & =\frac{1}{2}\label{eq:SRA_stab2_conds}\\
\alpha_{1}+\alpha_{2} & =1\nonumber \\
0<\left(A_{22}-A_{12}\right)\alpha_{1}+\left(A_{11}-A_{21}\right)\alpha_{2} & <1\nonumber
\end{align}
However, A-stability is not sufficient for most ODE integrators to
properly handle stiff equations and thus extra properties generally
imposed \cite{RN3790}. One important property we wish to extend to
stochastic integrators is L-stability. The straightforward extension
of L-stability is the condition
\begin{equation}
\lim_{z\rightarrow\infty}G(z)=0.\label{eq:SRA_L}
\end{equation}
This implies that
\begin{equation}
\frac{-A_{11}A_{22}+A_{11}\alpha_{2}+A_{12}A_{21}-A_{12}\alpha_{1}-A_{21}\alpha_{2}+A_{22}\alpha_{2}\alpha_{1}}{A_{12}A_{21}-A_{11}A_{22}}=0\label{eq:SRA_Lcond}
\end{equation}
The denominator is $-\det(A^{(0)})$ which implies $A^{(0)}$ must
be non-singular. Next, we attempt to impose B-stability on the drift
portion of the method. We use the condition due to Burrage and Butcher
that for $B=\text{diag}\left(\alpha_{1},\alpha_{2}\right)$ $M=BA^{(0)}+A^{(0)}B-\alpha\alpha^{T}$
(for ODEs) \cite{RN3791}, we require both $B$ and $M$ to be non-negative
definite. However, in the supplemental Mathematica notebooks we show
computationally that there is no 2-stage SRK method of this form which
satisfies all three of these stability conditions. Thus we settle
for A-stability and L-stability.

Recalling that $c^{(0)}$ and $c^{(1)}$ are the locations in time
where $f$ and $g$ are approximated respectively, we wish to impose
\begin{align}
c_{1}^{(0)} & =0\label{eq:SRA_ccond}\\
c_{2}^{(0)} & =1\nonumber \\
c_{1}^{(1)} & =0\nonumber \\
c_{2}^{(1)} & =1\nonumber
\end{align}
so that the error estimator covers the entire interval of integration.
Since $c^{(0)}=A^{(0)}e$, this leads to the condition $A_{21}+A_{22}=1$.
Using the constraint-satisfaction algorithm FindInstance in Mathematica,
we look for tableaus which satisfy the previous conditions with the
added constraint of semi-implicitness, i.e. $B^{(0)}$ is lower triangular.
This assumption is added because the inverse of the normal distribution
has unbounded moments, and thus in many cases it mathematically simpler
to consider the diffusion term as explicit (though there are recent
methods which drop this requirement via truncation or extra assumptions
on the solution \cite{RN3806}). However, we find that there is no
coefficient set which meets all of these requirements. However, if
we relax the interval estimate condition to allow $0\leq c_{2}^{(0)}\leq1$,
we find an A-L stable method:
\begin{align}
A^{(0)} & =\left(\begin{array}{cc}
1 & \frac{-41}{64}\\
\frac{32}{41} & \frac{9}{41}
\end{array}\right),\thinspace\thinspace\thinspace B^{(0)}=\left(\begin{array}{cc}
\frac{5}{8} & 0\\
0 & \frac{7}{3}
\end{array}\right),\thinspace\thinspace\thinspace\thinspace\alpha=\left(\begin{array}{c}
\frac{32}{41}\\
\frac{9}{41}
\end{array}\right)\label{eq:LSRA}\\
\beta^{(1)} & =\left(\begin{array}{c}
0\\
1
\end{array}\right),\thinspace\thinspace\thinspace\beta^{(2)}=\left(\begin{array}{c}
1\\
-1
\end{array}\right),\thinspace\thinspace\thinspace c^{(0)}=\left(\begin{array}{c}
\frac{23}{64}\\
1
\end{array}\right),\,\,\,c^{(1)}=\left(\begin{array}{c}
0\\
1
\end{array}\right)\nonumber
\end{align}
which we denote LSRA. If we attempt to look for a 2-stage SDIRK-like
method to reduce the complexity of the implicit equation, i.e. $A_{12}^{(0)}=0$,
using FindInstance we find the constraints unsatisfiable. Note that
if we drop the semi-implicit assumption we find that the full constraints
cannot be satisfied there (we still cannot satisfy $c_{1}^{(0)}=0$
and $c_{2}^{(0)}=1$), and there does not exist a 2-stage A-L stable
SDIRK method in that case.
\begin{center}
	\begin{figure}
		\begin{centering}
			\includegraphics[scale=0.5]{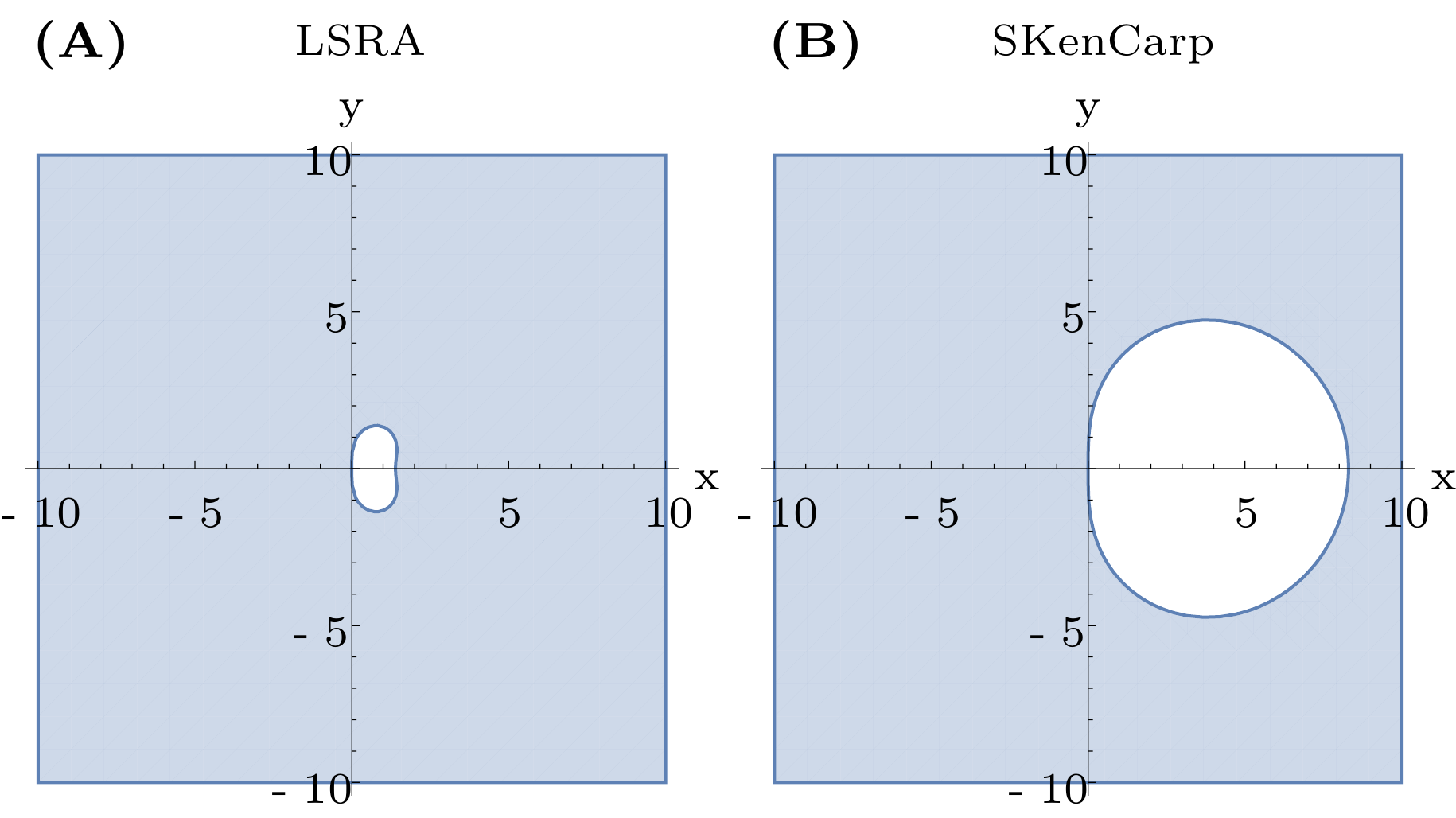}
			\par\end{centering}
		\caption{\textbf{Implicit SRA Stability Regions}. The stability regions ($\left|G(z)\right|<1$)
			are plotted in the $(x,y)$-plane for $z=x+iy$. \textbf{(A)} LSRA.
			\textbf{(B)} SKenCarp\label{fig:Implicit-Stability-Regions}}
	\end{figure}
	\par\end{center}

\subsubsection{Extensions of ODE Implicit Runge-Kutta Methods to Implicit SRA Methods}

Since the stability region of the SRA methods is completely determined
by the deterministic portion $A^{(0)}$, in some cases there may exist
a sensible extension of implicit Runge-Kutta methods for ordinary
differential equations to high order adaptive methods stochastic differential
equations with additive noise which keep the same stability properties.
Since the order constraints which only involve the deterministic portions
$A^{(0)}$, $c^{(0)}$, and $\alpha$ match the conditions required
for ODE integrators, existence is dependent on finding $\beta^{(1)}$,
$\beta^{(2)}$, $c^{(1)}$, and $B^{(0)}$ that satisfy the full order
constraints. In this case, an adaptive error estimator can be added
by using the same estimator as the ODE method (which we call $E_{D}$)
but adding the absolute size of the stochastic portions
\begin{equation}
E_{N}=\left|\sum_{i=1}^{s}\left(\beta_{i}^{(1)}I_{(1)}+\beta_{i}^{(2)}\frac{I_{(1,0)}}{h}\right)\right|\label{eq:EN}
\end{equation}
leading to the error estimator
\begin{equation}
E=\delta E_{D}+E_{N}.\label{eq:E}
\end{equation}
This can be shown similarly to the construction in \cite{RN3787}.
Given the large literature on implicit RK methods for ODEs, this presents
a large pool of possibly good methods and heuristically one may believe
that these would do very well in the case of small noise.

However, we note that there does not always exist such an extension.
Using the constraint-satisfaction algorithm FindInstance in Mathematica,
we looked for extensions of the explicit first stage singly-diagonally
implicit RK (ESDIRK) method TRBDF2 \cite{RN3793} and could not find
values satisfying the constraints. In addition, we could not find
values for an extension of the 5th order Radau IIA method \cite{RN3794,RN3790}
which satisfies the constraints. In fact, our computational search
could not find any extension of a 3-stage L-stable implicit RK method
which satisfies the constraints.

But, the 4-stage 3rd order ODE method due to Kennedy and Carpenter
\cite{RN3792} can be extended to the following:

\begin{align}
A^{(0)} & =\left(\begin{array}{cccc}
0 & 0 & 0 & 0\\
\frac{1767732205903}{4055673282236} & \frac{1767732205903}{4055673282236} & 0 & 0\\
\frac{2746238789719}{10658868560708} & -\frac{640167445237}{6845629431997} & \frac{1767732205903}{4055673282236} & 0\\
\frac{1471266399579}{7840856788654} & -\frac{4482444167858}{7529755066697} & \frac{11266239266428}{11593286722821} & \frac{1767732205903}{4055673282236}
\end{array}\right),\label{eq:KenCarp4}\\
\alpha & =\left(\begin{array}{c}
\frac{1471266399579}{7840856788654}\\
-\frac{4482444167858}{7529755066697}\\
\frac{11266239266428}{11593286722821}\\
\frac{1767732205903}{4055673282236}
\end{array}\right)\nonumber \\
\beta^{(1)} & =\left(\begin{array}{c}
0\\
0\\
0\\
1
\end{array}\right),\thinspace\thinspace\thinspace\beta^{(2)}=\left(\begin{array}{c}
1\\
0\\
0\\
-1
\end{array}\right),\thinspace\thinspace\thinspace c^{(0)}=\left(\begin{array}{c}
0\\
\frac{1767732205903}{4055673282236}\\
\frac{3}{5}\\
1
\end{array}\right),\,\,\,c^{(1)}=\left(\begin{array}{c}
0\\
0\\
0\\
1
\end{array}\right)\nonumber
\end{align}

{\scriptsize
\begin{align*}
B_{2,1}^{(0)} & \approx-12.246764387585055918338744103409192607986567514699471403397969732723452087723101\\
B_{4,3}^{(0)} & \approx-14.432096958608752822047165680776748797565142459789556194474191884258734697161106
\end{align*}
}%

The exact values for $B_{2,1}$ and $B_{4,3}$ are shown in \ref*{subsec:SKenCarpExact}.
(E)SDIRK methods are particularly interesting because these methods
can be solved using a single factorization of the function of the
Jacobian $I-\gamma dtJ$ where $J$ is the Jacobian. Additionally,
explicit handling of the noise term is similar to the Implicit-Explicit
(IMEX) form for additive Runge-Kutta methods in that it occurs by
adding a single constant term to the Newton iterations in each stage,
meaning it does not significantly increase the computational cost.
The chosen ESDIRK method has a complimentary explicit tableau to form
an IMEX additive Runge-Kutta method, and the chosen values for the
stochastic portions are simultaneously compatible with the order conditions
for this tableau. In Section \ref{subsec:Convergence-Tests} we numerically
investigate the order of the IMEX extension of the method and show
that it matches the convergence of the other SRA methods on the test
equation. One thing to note is that since the problem is additive
noise the method is never implicit in the dependent variables in the
noise part, so in theory this can also be extended with $B^{(1)}$
implicit as well (with convergence concerns due to the non-finite
inverse moments of the Normal distribution \cite{RN3169}).

\subsubsection{Note on Implementation and Stage Prediction}

One has to be careful with the implementation to avoid accumulation
of floating point error for highly stiff equations. For our implementation,
we used a method similar to that described in \cite{RN3793}. The
implicit stages were defined in terms of
\begin{equation}
z_{i}=hf\left(t+c_{i}^{(0)},H_{j}^{(0)}\right)\label{eq:zi}
\end{equation}
where $X_{0}$ is the previous step, and thus the iterations become
\begin{equation}
H_{j}^{(0)}=\gamma z_{i}+\sum_{j=1}^{i-1}A_{ij}^{(0)}f\left(t_{n}+c_{j}^{(0)}h,H_{j}^{(0)}\right)h+\sum_{j=1}^{i-1}B_{ij}^{(0)}g\left(t_{n}+c_{j}^{(1)}h\right)\frac{I_{(1,0)}}{h}=\gamma z_{i}+\alpha_{i}.\label{eq:z}
\end{equation}
This gives the implicit system for the residual:
\[
G(z_{i})=z_{i}-hf\left(t_{n}+c_{i}^{(0)}h,\gamma z_{i}+\alpha_{i}\right)
\]
which has a Jacobian $I-\gamma hJ$ where $J$ is the Jacobian of
$f$ and thus is the same for each stage. For choosing the values
to start the Newton iterations, also known as stage prediction, we
tested two methods. The first is the trivial stage predictor which
is $z_{i}=z_{j}$ for the $j$ s.t. $j<i$ and $c_{j}<c_{i}$, i.e.
using the closest derivative estimate. The other method that was tested
is what we denote the stochastic minimal residual estimate given by
$H_{j}^{(0)}=\alpha_{i}$ or $z_{i}=0$. This method takes into account
the stochastic bursts at a given step and thus demonstrated much better
stability.

\subsubsection{Note on Mass Matrices}

We note that these methods also apply to solving ODEs with mass-matrices
of the form:
\[
MdX_{t}=f(t,X_{t})dt+Mg(t,X_{t})dW_{t}.
\]
The derivation of the method is the same, except in this case we
receive the implicit system
\[
G(z_{i})=Mz_{i}-hf\left(t_{n}+c_{i}^{(0)}h,\gamma z_{i}+\alpha_{i}\right)
\]
which has a Jacobian $M-\gamma hJ$. Like in the ODE case, these
implicit methods can thus solve DAEs in mass-matrix form (the case
where $M$ is singular), though we leave discussion of convergence
for future research. One interesting property to note is that a zero
row in the mass matrix corresponds to a constraint equation which
is only dependent on the output of $f$ since the multiplication of
$g$ by $M$ is zero in that same corresponding row. Thus when a singular
mass matrix is applied to the noise equation, the corresponding constraints
are purely deterministic relations. Thus while this is a constrained
form, properties like conservation of energy in physical models can
still be placed on the solution using this mass-matrix formulation.

\section{Optimized-Stability Methods for Affine Noise via Transformation \label{sec:Affine-Transformation}}

Given the efficiency of the methods for additive noise, one method
for developing efficient methods for more general noise processes
is to use a transform of diagonal noise processes to additive noise.
This transform is due to Lamperti \cite{RN3525}, which states that
the SDE of the form

\begin{equation}
dX_{t}=f(t,X_{t})dt+\sigma(t,X_{t})R(t)dW_{t}\label{eq:Lamperti}
\end{equation}
where $\sigma>0$ is a diagonal matrix with diagonal elements $\sigma_{i}(t,X_{i,t})$
has the transformation
\begin{equation}
Z_{i,t}=\psi_{i}(t,X_{i,t})=\int\frac{1}{\sigma_{i}(x,t)}dx\mid_{x=X_{i,t}}\label{eq:LampertiIntegral}
\end{equation}
which will result in an Ito process with the $i$th element given
by
\begin{align}
dZ_{i,t}&=\left(\frac{\partial}{\partial t}\psi_{i}(t,x)\mid_{x=\psi^{-1}(t,Z_{i,t})}+\frac{f_{i}(t,\psi^{-1}(t,Z_{t}))}{\sigma_{i}\left(t,\psi_{i}^{-1}\left(t,Z_{i,t}\right)\right)}-\frac{1}{2}\frac{\partial}{\partial x}\sigma_{i}\left(t,\psi_{i}^{-1}\left(t,Z_{i,t}\right)\right)\right)dt\\
&+\sum_{j=1}^{n}r_{ij}(t)dw_{j,t}\label{eq:LampertiTransform}
\end{align}
with
\begin{equation}
X_{t}=\psi^{-1}\left(t,Z_{t}\right).\label{eq:LampertiInv}
\end{equation}
This is easily verified using Ito's Lemma. In the case of mixed multiplicative
and additive noise (affine noise), the vector equation:
\begin{equation}
dX_{t}=f(t,X_{t})dt+\left(\sigma_{M}X_{t}+\sigma_{A}\right)dW_{t}\label{eq:LampertiAffine}
\end{equation}
with $\sigma_{M}>0$ and $\sigma_{A}>0$, the transform becomes element-wise
in the system. Thus we can consider the one-dimensional case. Since
$\psi(t,X_{t})=\int\left(\frac{1}{\sigma_{M}X_{t}+\sigma_{A}}\right)dx\mid_{x=X_{t}}=\frac{\log(\sigma_{M}X_{t}+\sigma_{A})}{\sigma_{M}}$,
then $X_{t}=\frac{\exp\left(\sigma_{M}Z_{t}\right)-\sigma_{A}}{\sigma_{M}}$
and
\begin{align}
dZ_{t} & =\tilde{f}(t,X_{t})dt+dW_{t}\label{eq:LampertiAffineSol}\\
\tilde{f}(t,X_{t}) & =\frac{f(t,X_{t})}{\sigma_{M}X_{t}+\sigma_{A}}-\frac{1}{2}\sigma_{M}\nonumber
\end{align}
provided $\sigma_{M}X_{t}$ is guaranteed to be sufficiently different
from $\sigma_{A}$ to not cause definitional issues. It is common
in biological models like chemical reaction networks that $X_{t}\geq0$,
in which case this is well-defined for any $\sigma_{A}>0$ when $\sigma_{M}>0$.

For numerical problem solving environments (PSEs), one can make use
of this transformation in two ways. Source transformations could transform
affine noise SDEs element-wise to solve for the vector $Z_{t}$ which
is the same as $X_{t}$ if $\sigma_{M}\neq0$ and is the transformed
$X_{t}$ otherwise (assuming parameters must be positive). When doing
so, references of $X_{i,t}$ must be changed into $\frac{\exp(\sigma_{M}Z_{i,t})-\sigma_{A}}{\sigma_{M}}$.
For example, the affine noise Lotka-Volterra SDE:

\begin{align*}
dx & =\left(ax-bxy\right)dt+\left(\sigma_{M}x+\sigma_{A}\right)dW_{t}^{1}\\
dy & =\left(-cy+dxy\right)dt+\sigma_{\tilde{A}}dW_{t}^{2}
\end{align*}
only has noise on the first term, so this transforms to
\begin{align*}
x & =\frac{\exp(\sigma_{M}z)-\sigma_{A}}{\sigma_{M}}\\
dz & =\left(\frac{ax-bxy}{\sigma_{M}x+\sigma_{A}}-\frac{1}{2}\sigma_{M}\right)dt+dW_{t}^{1}\\
dy & =\left(-cy+dxy\right)dt+\sigma_{\tilde{A}}dW_{t}^{2}
\end{align*}
along with the change to the initial condition and can thus be solved
with the SRA methods. We note a word of caution that the above transformation
only holds when $\sigma_{A}>0$ and when $\sigma_{A}=0$, the transformation
is different, with $X_{t}=\frac{\exp\left(Z_{t}\right)}{\sigma_{M}}$
(instead of $\frac{\exp\left(\sigma_{M}Z_{t}\right)}{\sigma_{M}}$
which one would get by taking $\sigma_{A}=0$).

Instead of performing the transformations directly on the functions
themselves, we can modify the SRA algorithm to handle this case as:

{\scriptsize
\begin{equation}
X_{n+1}=\psi^{-1}\left(\psi\left(X_{n}\right)+\sum_{i=1}^{s}\alpha_{i}\tilde{f}\left(t_{n}+c_{i}^{(0)}h,H_{i}^{(0)}\right)+\sum_{i=1}^{s}\left(\beta_{i}^{(1)}I_{(1)}+\beta_{i}^{(2)}\frac{I_{(1,0)}}{h}\right)\tilde{g}\left(t_{n}+c_{i}^{(1)}h\right)\right)\label{eq:AffineStep}
\end{equation}
}%

with stages
\begin{equation}
H_{i}^{(0)}=\psi\left(X_{n}\right)+\sum_{j=1}^{s}A_{ij}^{(0)}\tilde{f}\left(t_{n}+c_{j}^{(0)}h,H_{j}^{(0)}\right)h+\sum_{j=1}^{s}B_{ij}^{(0)}\tilde{g}\left(t_{n}+c_{j}^{(1)}h\right)\frac{I_{(1,0)}}{h}\label{eq:AffineStep2}
\end{equation}
where $\psi$ is the element-wise function:
\[
\psi_{i}(x)=\begin{cases}
\frac{\log(\sigma_{i,M}x+\sigma_{i,A})}{\sigma_{i,M}} & \sigma_{i,M}>0,\sigma_{i,A}>0\\
\frac{\log(x)}{\sigma_{i,M}} & \sigma_{i,M}>0,\sigma_{i,A}=0\\
x & o.w.
\end{cases},
\]
\[
\thinspace\thinspace\thinspace\psi_{i}^{-1}(z)=\begin{cases}
\frac{\exp\left(\sigma_{i,M}z\right)}{\sigma_{i,M}} & \sigma_{i,M}>0,\sigma_{i,A}>0\\
\frac{\exp\left(z\right)}{\sigma_{i,M}} & \sigma_{i,M}>0,\sigma_{i,A}=0\\
x & o.w.
\end{cases}
\]
and
\[
\tilde{g}_{i}(t)=\begin{cases}
1 & \sigma_{i,M}>0\\
\sigma_{i,A} & o.w.
\end{cases}
\]
This can be summarized as performing all internal operations in $Z$-space
(where the equation is additive) but saving each step in $X$-space.

\section{Optimized-Stability Order 1.5 SRK Methods with Diagonal Noise}

\subsection{The Stability Equation for Order 1.5 SRK Methods with Diagonal Noise}

For diagonal noise, we use the mean-square definition of stability
\cite{RN3169}. A method is mean-square stable if $\lim_{n\rightarrow\infty}\mathbb{E}\left(\left|X_{n}\right|^{2}\right)=0$
on the test equation
\begin{equation}
dX_{t}=\mu X_{t}dt+\sigma X_{t}dW_{t}.\label{eq:multtest}
\end{equation}
In matrix form we can re-write our method as given by

\begin{align}
X_{n+1}&=X_{n}+\mu h\left(\alpha\cdot H^{(0)}\right)+\sigma I_{(1)}\left(\beta^{(1)}\cdot H^{(1)}\right)+\sigma\frac{I_{(1,1)}}{\sqrt{h}}\left(\beta^{(2)}\cdot H^{(1)}\right)\\
&+\sigma\frac{I_{(1,0)}}{h}\left(\beta^{(3)}\cdot H^{(1)}\right)+\sigma\frac{I_{(1,1,1)}}{h}\left(\beta^{(4)}\cdot H^{(1)}\right)\label{eq:matrix_method}
\end{align}
with stages

\begin{eqnarray}
H^{(0)} & = & X_{n}+\mu\Delta tA^{(0)}H^{(0)}+\sigma\frac{I_{(1,0)}}{h}B^{(0)}H^{(1)},\label{eq:matrix_stages}\\
H^{(1)} & = & X_{n}+\mu\Delta tA^{(1)}H^{(0)}+\sigma\sqrt{\Delta t}B^{(1)}H^{(1)}\nonumber
\end{eqnarray}
where $\hat{X_{n}}$ is the size $s$ constant vector of $X_{n}$.

\begin{eqnarray}
H^{(0)} & = & \left(I-hA^{(0)}\right)^{-1}\left(\hat{X_{n}}+\sigma\frac{I_{(1,0)}}{h}B^{(0)}H^{(1)}\right),\label{matrx_stages_solve}\\
H^{(1)} & = & \left(I-\sigma\sqrt{h}B^{(1)}\right)^{-1}\left(\hat{X_{n}}+\mu hA^{(1)}H^{(0)}\right)\nonumber
\end{eqnarray}
By the derivation in the appendix, we receive the equation

{\tiny{}
	\begin{eqnarray}
	S=E\left[\frac{U_{n+1}^{2}}{U_{n}^{2}}\right] & = & \{1+\mu ht\left(\alpha\cdot\left[\left(I-\mu\Delta tA^{(0)}-\mu\sigma I_{(1,0)}A^{(1)}B^{(0)}\left(I-\sigma\sqrt{h}B^{(1)}\right)^{-1}\right)^{-1}\left(I+\sigma\frac{I_{(1,0)}}{h}B^{(0)}\left(I-\sigma\sqrt{h}B^{(1)}\right)^{-1}\right)\right]\right)\label{eq:stability_solve}\\
	&  & +\sigma I_{(1)}\left(\beta^{(1)}\cdot\left[\left(I-\sigma\sqrt{h}B^{(1)}-\mu hA^{(1)}\left(I-\mu hA^{(0)}\right)^{-1}\sigma\frac{I_{(1,0)}}{h}B^{(0)}\right)^{-1}\left(I+\mu hA^{(1)}\left(I-\mu hA^{(0)}\right)^{-1}\right)\right]\right)\nonumber \\
	&  & +\sigma\frac{I_{(1,1)}}{\sqrt{h}}\left(\beta^{(2)}\cdot\left[\left(I-\sigma\sqrt{h}B^{(1)}-\mu hA^{(1)}\left(I-\mu hA^{(0)}\right)^{-1}\sigma\frac{I_{(1,0)}}{h}B^{(0)}\right)^{-1}\left(I+\mu hA^{(1)}\left(I-\mu hA^{(0)}\right)^{-1}\right)\right]\right)\nonumber \\
	&  & +\sigma\frac{I_{(1,0)}}{h}\left(\beta^{(3)}\cdot\left[\left(I-\sigma\sqrt{h}B^{(1)}-\mu hA^{(1)}\left(I-\mu hA^{(0)}\right)^{-1}\sigma\frac{I_{(1,0)}}{h}B^{(0)}\right)^{-1}\left(I+\mu hA^{(1)}\left(I-\mu hA^{(0)}\right)^{-1}\right)\right]\right)\nonumber \\
	&  & +\sigma\frac{I_{(1,1,1)}}{h}\left(\beta^{(4)}\cdot\left[\left(I-\sigma\sqrt{h}B^{(1)}-\mu hA^{(1)}\left(I-\mu hA^{(0)}\right)^{-1}\sigma\frac{I_{(1,0)}}{h}B^{(0)}\right)^{-1}\left(I+\mu hA^{(1)}\left(I-\mu hA^{(0)}\right)^{-1}\right)\right]\right)\}^{2}\nonumber
	\end{eqnarray}
}We apply the substitutions from the Appendix and let

\begin{align}
z & =\mu h,\label{eq:axis}\\
w & =\sigma\sqrt{h}.\nonumber
\end{align}
In this space, $z$ is the stability variable for the drift term and
$w$ is the stability in the diffusion term. Under this scaling $\left(h,\sqrt{h}\right)$,
the equation becomes independent of $h$ and thus becomes a function
$S(z,w)$ on the coefficients of the SRK method where mean-square
stability is achieved when $\left|S(z,w)\right|<1$. The equation
$S(z,w)$ in terms of its coefficients for explicit methods ($A^{(i)}$
and $B^{(i)}$ lower diagonal) has millions of terms and is shown
in the supplemental Mathematica notebook. Determination of the stability
equation for the implicit methods was found to be computationally
intractable and is an avenue for further research.

\subsection{An Optimization Problem for Determination of Coefficients}

We wish to determine the coefficients for the diagonal SRK methods
which optimize the stability. To do so, we generate an optimization
problem which we can numerically solve for the coefficients. To simplify
the problem, we let $z,w\in\mathbb{R}$. Define the function

\begin{equation}
f\left(z,w;N,M\right)=\int_{-M}^{M}\int_{-N}^{1}\chi_{S(z,w)\leq1}(z,w)dzdw.\label{eq:SRIarea}
\end{equation}
Notice that for $N,M\rightarrow\infty$, $f$ is the area of the stability
region. Thus we define the stability-optimized diagonal SRK method
as the set of coefficients which achieves

\begin{align}
\max_{A^{(i)},B^{(i)},\beta^{(i)},\alpha} & f(z,w)\label{eq:opt}\\
\text{subject to: } & \text{Order Constraints}\nonumber
\end{align}
However, like with the SRK methods for additive noise, we impose
a few extra constraints to add robustness to the error estimator.
In all cases we impose $0<c_{i}^{(0)},c_{i}^{(1)}<1$ . Additionally
we can prescribe $c_{4}^{(0)}=c_{4}^{(1)}=1$ which we call the End-C
Constraint. Lastly, we can prescribe the ordering constraint $c_{1}^{(j)}<c_{2}^{(j)}<c_{3}^{(j)}<c_{4}^{(j)}$
which we denote as the Inequality-C Constraint.

The resulting problem is a nonlinear programming problem with 44 variables
and 42-48 constraint equations. The objective function is the two-dimensional
integral of a discontinuous function which is determined by a polynomial
of in $z$ and $w$ with approximately 3 million coefficients. To
numerically approximate this function, we calculated the characteristic
function on a grid with even spacing $dx$ using a CUDA kernel and
found numerical solutions to the optimization problem using the JuMP
framework \cite{RN3795} with the NLopt backend \cite{RN3788}. A
mixed approach using many solutions of the semi-local optimizer LN\_AUGLAG\_EQ
\cite{RN3797,RN3798} and fewer solutions from the global optimizer
GN\_ISRES \cite{RN3796} were used to approximate the optimality of
solutions. The optimization was run many times in parallel until many
results produced methods with similar optimality, indicating that
we likely obtained values near the true minimum.

The parameters $N$ and $M$ are the bounds on the stability region
and also represent a trade-off between the stability in the drift
and the stability in the diffusion. A method which is optimized when
$M$ is small would be highly stable in the case of small noise, but
would not be guaranteed to have good stability properties in the presence
of large noise. Thus these parameters are knobs for tuning the algorithms
for specific situations, and thus we solved the problem for different
combinations of $N$ and $M$ to determine different algorithms for
the different cases.

\subsection{Resulting Approximately-Optimal Methods}

The coefficients generated for approximately-optimal methods fall
into three categories. In one category we have the drift-dominated
stability methods where large $N$ and small $M$ was optimized. On
the other end we have the diffusion-dominated stability methods where
large $M$ and small $N$ was optimized. Then we have the mixed stability
methods which used some mixed size choices for $N$ and $M$. As a
baseline, we optimized the objective without constraints on the $c_{i}$
to see what the ``best possible method'' would be. When this was
done with large $N$ and $M$, the resulting method, which we name
SOSRI, has almost every value of $c$ satisfy the constraints, but
with $c_{2}^{(0)}\approx-0.04$ and $c_{4}^{(0)}\approx3.75$. To
see if we could produce methods which were more diffusion-stable,
we decreased $N$ to optimize more in $w$ but failed to produce methods
with substantially enlarged diffusion-stability over SOSRI.

Adding only the inequality constraints on the $c_{i}$ and looking
for methods for drift-dominated stability, we failed to produce methods
whose $c_{i}$ estimators adequately covered the interval. Some of
the results did produce stability regions similar to SOSRI but with
$c_{i}^{(0)}<0.5$ which indicates the method could have problems
with error estimation. When placing the equality constraints on the
edge $c_{i}$, one method, which we label SOSRI2, resulted in similar
stability to SOSRI but satisfy the $c_{i}$ constraints. In addition,
this method satisfies $c_{3}^{(0)}=c_{4}^{(0)}=1$ and $c_{3}^{(1)}=c_{4}^{(1)}=1$,
a property whose use will be explained in Section \ref{subsec:StabilityDetection}.
The stability regions for these methods is shown in Figure \ref{fig:SOSRI-Stability-Regions.}.

To look for more diffusion-stable methods, we dropped to $N=6$ to
encourage the methods to expand the stability in the $w$-plane. However,
we could not find a method whose stability region went substantially
beyond $\left[-2,2\right]$ in $w$. This was further decreased to
$N=1$ where methods still could not go substantially beyond $\left|2\right|$.
Thus we were not able to obtain methods optimized for the diffusion-dominated
case. This hard barrier was hit under many different constraint and
objective setups and under thousands of optimization runs, indicating
there might be a diffusion-stability barrier for explicit methods.
\begin{center}
	\begin{figure}
		\begin{centering}
			\includegraphics[scale=0.46]{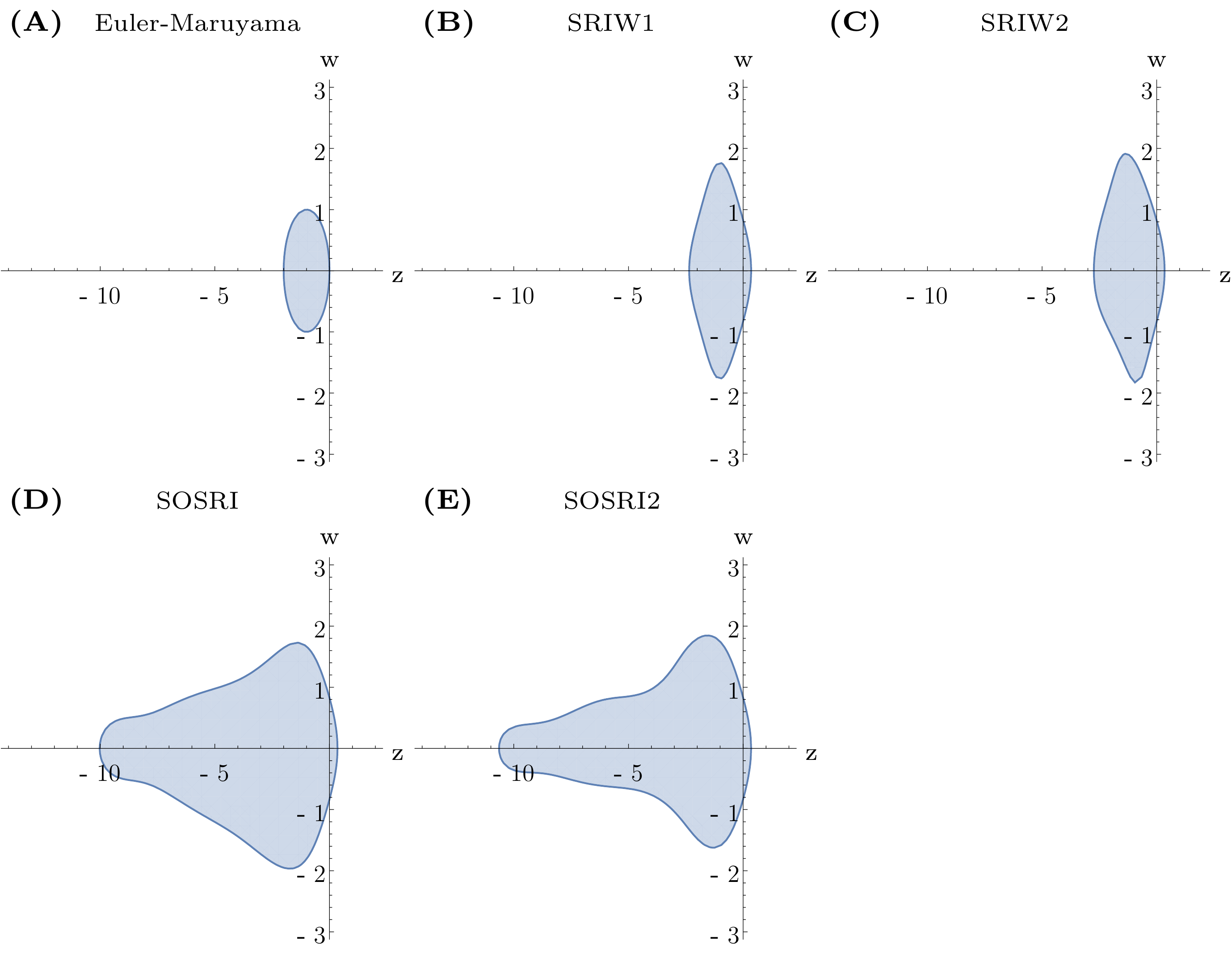}
			\par\end{centering}
		\caption{\textbf{SOSRI Stability Regions}. The stability regions $\left(S(z,w)\leq1\right)$for
			the previous and SOSRI methods are plotted in the $(z,w)$-plane.
			\textbf{(A) }Euler-Maruyama. \textbf{(B) }SRIW1. \textbf{(C) }SRIW2.
			\textbf{(D) }SOSRI. \textbf{(E) }SOSRI2\label{fig:SOSRI-Stability-Regions.}}
	\end{figure}
	\par\end{center}

\subsection{Approximately-Optimal Methods with Stability Detection and \\
	Switching Behaviors\label{subsec:StabilityDetection}}

In many real-world cases, one may not be able to clearly identify
a model as drift-stability bound or diffusion-stability bound, or
if the equation is stiff or non-stiff. In fact, many models may switch
between such extremes. An example is a model with stochastic switching
between different steady states. In this case, we have that the diffusion
term $f(t,X_{ss})\approx0$ in the area of many stochastic steady
states, meaning that while straddling a steady state the integration
is heavily diffusion-stability dominated and usually non-stiff. However,
when switching between steady states, $f$ can be very large and stiff,
causing the integration to be heavily drift-stability dominated. Since
these switches are random, the ability to adapt between these two
behaviors could be key to achieving optimal performance. Given the
trade-off, we investigated how our methods allow for switching between
methods which optimize for the different situations.

The basis for our method is an extension of a method proposed for
deterministic differential equations \cite{RN3527,RN3526,RN3790}.
The idea is to create a cheap approximation to the dominant eigenvalues
of the Jacobians for the drift and diffusion terms. If $v$ is the
eigenvector of the respective Jacobian, then for $\Vert v\Vert$ sufficiently
small,
\begin{align}
\left|\lambda_{D}\right|&\approx\frac{\Vert f(t,x+v)-f(t,x)\Vert}{\Vert v\Vert},\\
\left|\lambda_{N}\right|&\approx\frac{\Vert g(t,x+v)-g(t,x)\Vert}{\Vert v\Vert}\label{eq:eigen}
\end{align}
where $\left|\lambda_{D}\right|$ and $\left|\lambda_{N}\right|$
are the estimates of the dominant eigenvalues for the deterministic
and noise functions respectively. We have in approximation that $H_{i}^{(k)}$
is an approximation for $X_{t+c_{i}^{(k)}h}$ and thus the difference
between two successive approximations at the same time-point, $c_{i}^{(k)}=c_{j}^{(k)}$,
then the following serves as a local Jacobian estimate:

\begin{align}
\left|\lambda_{D}\right|&\approx\frac{\Vert f(t+c_{i}^{(0)}h,H_{i}^{(0)})-f(t+c_{j}^{(0)}h,H_{j}^{(0)})\Vert}{\Vert H_{i}^{(0)}-H_{j}^{(0)}\Vert},\\
\left|\lambda_{N}\right|&\approx\frac{\Vert f(t+c_{i}^{(1)}h,H_{i}^{(1)})-f(t+c_{j}^{(1)}h,H_{j}^{(1)})\Vert}{\Vert H_{i}^{(1)}-H_{j}^{(1)}\Vert}\label{eq:local_eigen}
\end{align}
If we had already computed a successful step, we would like to know
if in the next calculation we should switch methods due to stability.
Thus it makes sense to approximate the Jacobian at the end of the
interval, meaning $i=s$ and $j=s-1$ where $s$ is the number of
stages. Then if $z_{min}$ is the minimum $z\in\mathbb{R}$ such that
$z$ is in the stability region for the method, $\frac{h\left|\lambda_{D}\right|}{z_{min}}>1$
when the steps are outside the stability region. Because the drift
and mixed stability methods do not track the noise axis directly,
we instead modify $w_{min}$ to be $\frac{2}{3}$ of the maximum of
the stability region in the noise axis.

Hairer noted that, for ODEs, if a RK method has $c_{i}=c_{j}=1$,
then it follows that
\begin{equation}
\rho=\frac{\Vert k_{i}-k_{j}\Vert}{\Vert g_{i}-g_{j}\Vert}\label{eq:rho}
\end{equation}
where $k_{i}=f(t+c_{i}h,g_{i})$ is an estimate of the eigenvalues
for the Jacobian of $f$. Given the construction of SOSRI2, a natural
extension is

\begin{align}
\left|\lambda_{D}\right|&\approx\frac{\Vert f\left(t_{n}+c_{4}^{(0)}h,H_{4}^{(0)}\right)-f\left(t_{n}+c_{3}^{(0)}h,H_{3}^{(0)}\right)\Vert}{\Vert H_{4}^{(0)}-H_{3}^{(0)}\Vert},\\
\left|\lambda_{N}\right|&\approx\frac{\Vert g\left(t_{n}+c_{4}^{(0)}h,H_{4}^{(1)}\right)-g\left(t_{n}+c_{3}^{(0)}h,H_{3}^{(1)}\right)\Vert}{\Vert H_{4}^{(1)}-H_{3}^{(1)}\Vert}\label{eq:sriopt2_eigen}
\end{align}
Given that these values are all part of the actual step calculations,
this stiffness estimate essentially is free. By comparing these values
to the stability plot in Figure \ref{fig:SOSRA-Stability-Regions.},
we use the following heuristic to decide if SOSRI2 is stability-bound
in its steps:
\begin{enumerate}
	\item If $10>\left|\lambda_{D}\right|>2.5$, then we check if $h\left|\lambda_{N}\right|>\omega$.
	\item If $\left|\lambda_{D}\right|<2.5$, then we check if $h\left|\lambda_{N}\right|/2>\omega$.
\end{enumerate}
The denominator is chosen as a reasonable box approximation to the
edge of the stability region. $\omega$ is a safety factor: in theory
$\omega$ is $1$ since we divided by the edge of the stability region,
but in practice this is only an eigenvalue estimate and thus $\omega$
allows for a trade-off between the false positive and false negative
rates. If either of those conditions are satisfied, then $h$ is constrained
by the stability region. The solver can thus alert the user that the
problem is stiff or use this estimate to switch to a method more suitable
for stiff equations. In addition, the error estimator gives separate
error estimates in the drift and diffusion terms. A scheme could combine
these two facts to develop a more robust stiffness detection method,
and label the stiffness as either drift or diffusion dominated.

We end by noting that SOSRA2 has the same property, allowing stiffness
detection via

\begin{equation}
\left|\lambda_{D}\right|\approx\frac{\Vert f\left(t_{n}+c_{3}^{(0)}h,H_{3}^{(0)}\right)-f\left(t_{n}+c_{2}^{(0)}h,H_{2}^{(0)}\right)\Vert}{\Vert H_{3}^{(0)}-H_{2}^{(0)}\Vert}\label{eq:SOSRA_eigen}
\end{equation}
and, employing a similar method as the deterministic case, check
for stiffness via the estimate $h\left|\lambda_{D}\right|/5>\omega$.

In addition, stiff solvers can measure the maximal eigenvalues directly
from the Jacobian. Here we suggest the measure from Shampine \cite{RN3527,RN3526,RN3790}
of using $\Vert J\Vert_{\infty}$ as a cheap upper bound. For semi-implicit
methods like LSRA we only get a stability bound on the drift term,
but this should be sufficient since for additive noise diffusive noise
instability is not an issue.

\section{Numerical Results}

\subsection{Convergence Tests \label{subsec:Convergence-Tests}}

In order to test the efficiency and correctness of the SRA algorithms,
we chose to use the additive noise test Equation \ref{eq:ex1}. Figure
\ref{fig:convergence}A demonstrates that the SOSRA and SKenCarp methods
achieve the strong order 2.0 on Equation \ref{eq:ex3}. To test the
convergence of the SRI algorithms, we used the linear test Equation
\ref{eq:ex1}. Figure \ref{fig:convergence}B demonstrates that the
SOSRI methods achieve the strong order 1.5 on Equation \ref{eq:ex1}.
Lastly, we tested the convergence of the IMEX version of the SKenCarp
integrator. We defined the split SDE \ref{eq:ex3-1} as a modification
of Equation \ref{eq:ex3} where the $f_{1}$ part is solved implicitly
and the $f_{2}$ part is solved explicitly. Figure \ref{fig:convergence}C
demonstrates that the IMEX SKenCarp method achieves strong order 2.0
. Note that this does not demonstrate that the method always achieves
strong order 1.5 since sufficient conditions for the IMEX pairing
are unknown, but it gives numerical evidence that the method can be
high order.
\begin{center}
	\begin{figure}
		\begin{centering}
			\includegraphics[scale=0.5]{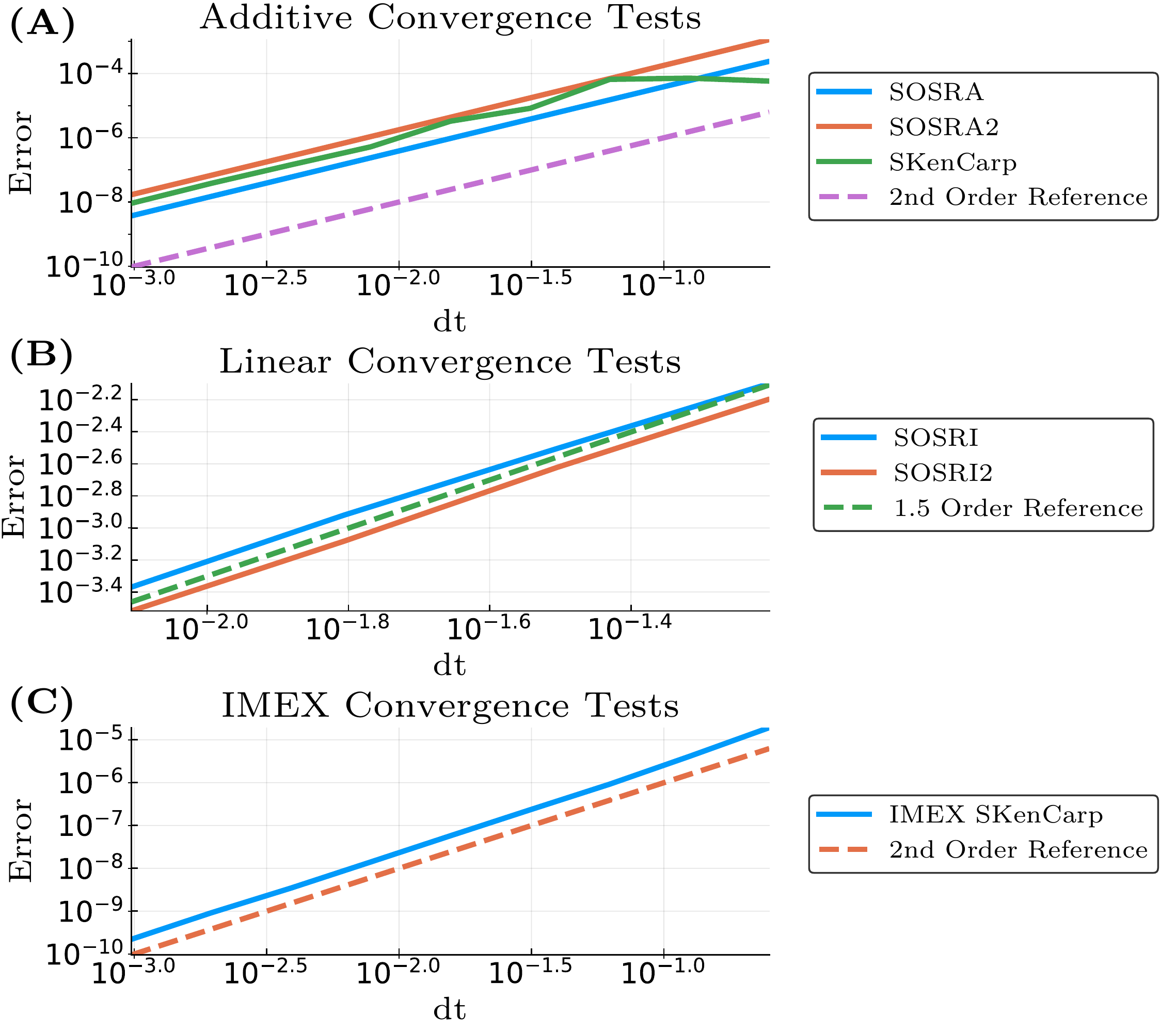}
			\par\end{centering}
		\caption{\textbf{Additive Noise Convergence Tests.} The error is averaged over
			1000 trajectories. Shown are the strong $l_{2}$ error along the time
			series of the solution.\textbf{ (A)} Convergence results on Equation
			\ref{eq:ex3}. The test used a fixed time step $h=1/2^{-2}$ to $h=1/2^{-10}$.\textbf{
				(B) }Convergence results on Equation \ref{eq:ex1}. The test used
			a fixed time step $h=1/2^{-4}$ to $h=1/2^{-7}$.\textbf{ (C)} Convergence
			results on the IMEX Equation \ref{eq:ex3-1}. The test used a fixed
			time step $h=1/2^{-2}$ to $h=1/2^{-10}$. \label{fig:convergence}}
	\end{figure}
	\par\end{center}

\subsection{SOSRA Numerical Efficiency Experiments}

\subsubsection{Additive Noise Lotka-Volterra (2 Non-Stiff SDEs)}

To test the efficiency we first plotted work-precision \cite{RN3353,RN3807,RN3790}
diagrams for the SOSRA, SOSRA2, and SKenCarp methods against the SRA1,
SRA2, SRA3 \cite{RN2707} methods, and fixed time step Euler-Maruyama
method (Milstein is equivalent to Euler-Maruyama in this case \cite{RN3169}).
We tested the error and timing on Equation \ref{eq:ex3}. In addition,
we tested using the Lotka-Volterra equation with additive noise Equation
\ref{eq:lotka_add}. Since \ref{eq:lotka_add} does not have an analytical
solution, a reference solution was computed using a low tolerance
solution via SOSRA for each Brownian trajectory. The plots show that
there is a minimal difference in efficiency between the SRA algorithms
for errors in the interval $\left[10^{-6},10^{-2}\right]$, while
these algorithms are all significantly more efficient than the Euler-Maruyama
method when the required error is $<10^{-4}$ (Figure \ref{fig:SOSRA-efficiency}).
The weak error work-precision diagrams show that when using between
100 to 10,000 trajectories, the weak error is less than the sample
error in the regime where there is no discernible efficiency difference
between the SRA methods. These results show that in the regime of
mild accuracy on non-stiff equations, the SOSRA, SOSRA2, and SKenCarp
methods are much more efficient than low order methods yet achieve
the same efficiency as the non-stability optimized SRA variants. Note
that these results also show that the error estimator for adaptivity
is highly conservative, generating solutions with around 2 orders
of magnitude less error than the tolerance suggests.
\begin{center}
	\begin{figure}
		\begin{centering}
			\includegraphics[scale=0.34]{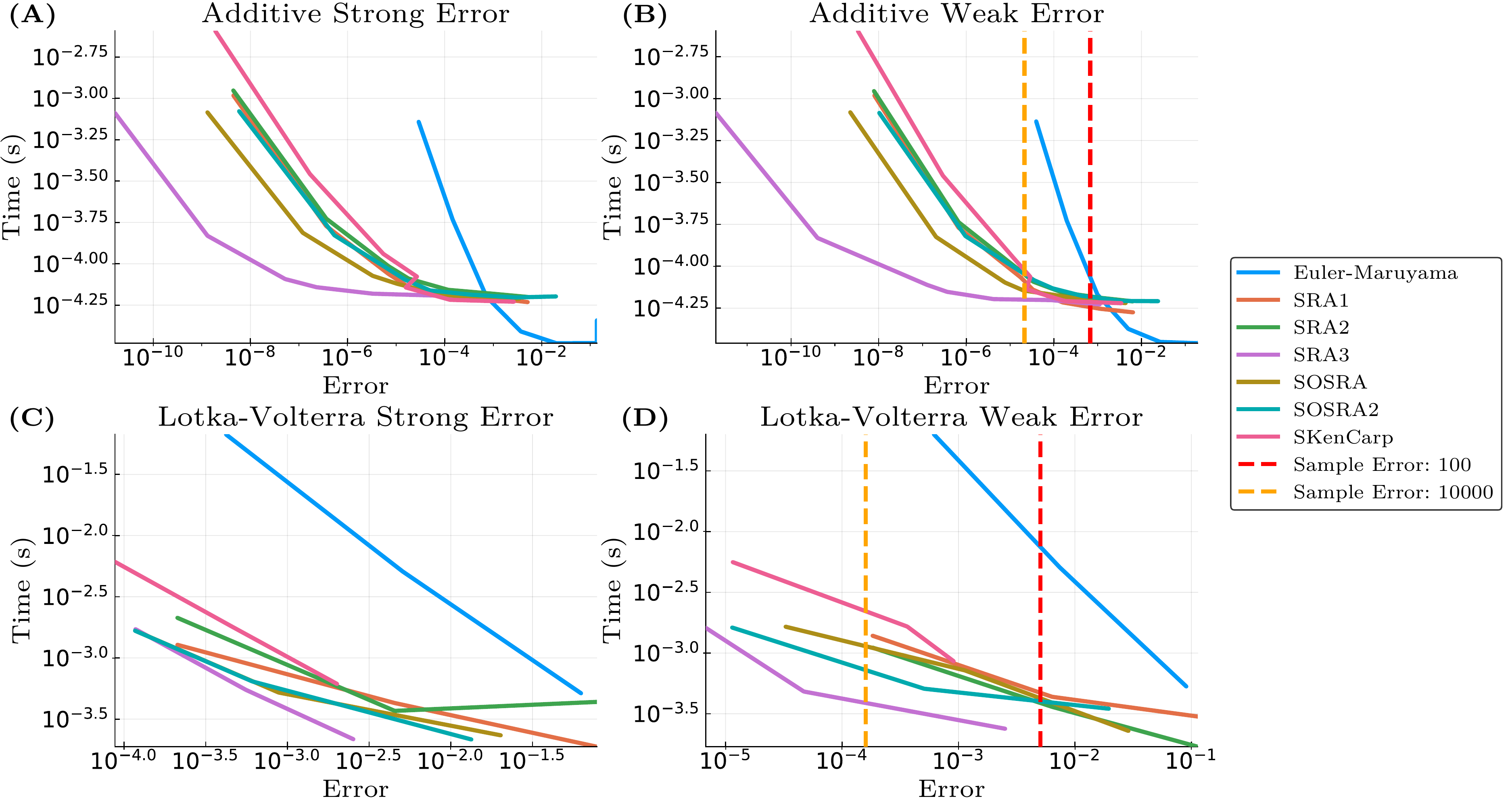}
			\par\end{centering}
		\caption{\textbf{SOSRA Efficiency Tests.} The error was taken as the average
			of 10,000 trajectories for Equation \ref{eq:ex3} and 100 trajectories
			for the Lokta-Volterra Equation \ref{eq:lotka_add}. The sample error
			was determined for the weak error as the normal 95\% confidence interval
			for the mean using the variance of the true solution Equation \ref{eq:sol3}
			or the variance of the estimated true solutions via low tolerance
			solutions. The time is the average time to compute a trajectory and
			is averaged over 1000 runs at the same tolerance or step size.\textbf{
				(A)} Shown are the work-precision plots for the methods on Equation
			\ref{eq:ex3}. Each of the adaptive time-stepping methods solved the
			problem on the interval using changing values of tolerances, with
			$tol=abstol=reltol$ starting at $10^{2}$ and ending at $10^{-4}$
			going in increments of $10$. The fixed time-stepping methods used
			time steps of size $h=1/5^{-1}$ to $h=1/5^{4}$, changing the value
			by factors of $5$. The error is the strong $l_{2}$ error computed
			over the time series. \textbf{(B)} Same setup as the previous plot
			but using the weak error at the final time-point. \textbf{(C) }Shown
			are the work-precision plots for the methods on the Equation \ref{eq:lotka_add}.
			Each of the adaptive time-stepping methods solved the problem on the
			interval using changing values of tolerances, with $tol=abstol=reltol$
			starting at $4^{-2}$ and ending at $4^{-4}$ going in increments
			of $4$. The fixed time-stepping methods used time steps of size $h=1/12^{-2.5}$
			to $h=1/12^{-6.5}$, changing the value by factors of $12$. The error
			is the strong $l_{2}$ error computed over the time series. \textbf{(D)}
			Same setup as the previous plot but using the weak error at the final
			time-point. \label{fig:SOSRA-efficiency}}
	\end{figure}
	\par\end{center}

\subsubsection{Addtive Noise Van Der Pol (2 Stiff SDEs)}

To test how efficiently the algorithms could achieve solve stiff equations,
we chose to analyze the qualitative results of the driven Van der
Pol equation. The driven Van der Pol equation is given by Equation
\ref{eq:van_der_pol} where $\mu$ is the driving factor. As $\mu$
increases the equation becomes more stiff. $\mu=10^{6}$ is a common
test for stiff ODE solvers \cite{RN3794}, with lower values used
to test the semi-stiff regime for ODEs. For our purposes, we chose
$\mu=10^{5}$ as a semi-stiff test case. The ODE case, solved using
the Tsit5 explicit Runge-Kutta algorithm \cite{RN3523,RN3784}, and
demonstrates the mild stiffness which is still well-handled by explicit
methods (Figure \ref{fig:AdditiveVanDerPol}A). We extend this model
to the driven Van der Pol model with additive noise Equation \ref{eq:van_der_pol_add}
where $\rho=3.0$ is the noise gain and $dW^{(1)}$ and $dW^{(2)}$
are independent Brownian motions. The solution to this model is interesting
because it gives the same qualitative behavior, large bursts when
$x(t)$ crosses zero, but in this case the zero crossings are stochastic.
Even at high tolerances, ($abstol=10$,$reltol=1/2^{1}$), SOSRA is
able to reproduce this qualitative behavior of the low tolerance solutions
(Figure \ref{fig:AdditiveVanDerPol}B), and SOSRA2 producing similar
results at the same tolerances a factor of two lower. Given the conservativeness
of the error estimators shown in previous (and other tests), this
case corresponds to roughly two decimal places of accuracy, which
is more than sufficient for many phenomenological models. However,
even at tolerances of $abstol=1/2^{3}$,$reltol=1/2^{3}$ SRA3 was
unable to reproduce the correct qualitative behavior (Figure \ref{fig:AdditiveVanDerPol}C).
Thus we decreased the tolerances by factors of 2 until it was able
to reproduce the correct qualitative results (Figure \ref{fig:AdditiveVanDerPol}D).
This shows that the SOSRA are more reliable on models with transient
stiffness. To test the impact on the run time of the algorithms, each
of the algorithms were run 100 times with the tolerance setup that
allows them to most efficiently generate correct qualitative results.
The run times are shown in Table \ref{tab:SRA-Runtimes-on-Van}, which
show that SRA1 takes more than 10 times and SRA3 nearly 4 times as
long as the SOSRA methods. In this case the implicit method SKenCarp
is the fastest by besting the SOSRA methods by more than 8x while
achieving similar qualitative results. This shows that as stiffness
comes into play, the SOSRA methods along with the implicit SKenCarp
method are more robust and efficient. The fixed time step methods
were far less efficient. Adaptive timestepping via rejection sampling
was crucial to the success of the SKenCarp method because it required
the ability to pull back to a smaller timestep when Newton iterations
diverged, otherwise it resulted in time estimates around 5x slower
than SOSRA.
\begin{center}
	\begin{table}
		\begin{centering}
			\begin{tabular}{|c|c|c|}
				\hline
				Algorithm & Run-time (seconds) & Relative Time (vs SKenCarp)\tabularnewline
				\hline
				\hline
				SKenCarp & 37.23 & 1.0x\tabularnewline
				\hline
				SOSRA & 315.58 & 8.5x\tabularnewline
				\hline
				SOSRA2 & 394.82 & 10.6x\tabularnewline
				\hline
				SRA3 & 1385.66 & 37.2x\tabularnewline
				\hline
				SRA1 & 3397.66 & 91.3x\tabularnewline
				\hline
				Euler-Maruyama & 5949.19 & 159.8x\tabularnewline
				\hline
				DISTM $\left(\theta=\frac{1}{2}\right)$ & 229111.15 & 6153x\tabularnewline
				\hline
			\end{tabular}
			\par\end{centering}
		\caption{\textbf{SRA Run times on Van der Pol with additive noise.} The additive
			noise Van der Pol equation was solved 100 times using the respective
			algorithms at the highest tolerance by powers of two which match the
			low tolerance solution to plotting accuracy. The fixed time step methods
			had their $\Delta t$ determined as the largest $\Delta t$ in increments
			of powers of 2 that produced no unstable trajectories. This resulted
			in $\Delta t=5e-8$ for both the Euler-Maruyama the Drift-Implicit
			Stochastic $\theta$-methods. Note that the total time of the drift-implicit
			stochastic $\theta$-method and the Euler-Maruyama method were determined
			by extrapolating the time from a single stable trajectory on $t\in[0,1]$
			due to time constraints. DISTM is the Drift-Implicit Stochastic
			$\theta$-Method\label{tab:SRA-Runtimes-on-Van}}
	\end{table}
	\par\end{center}

\begin{center}
	\begin{figure}
		\begin{centering}
			\includegraphics[scale=0.42]{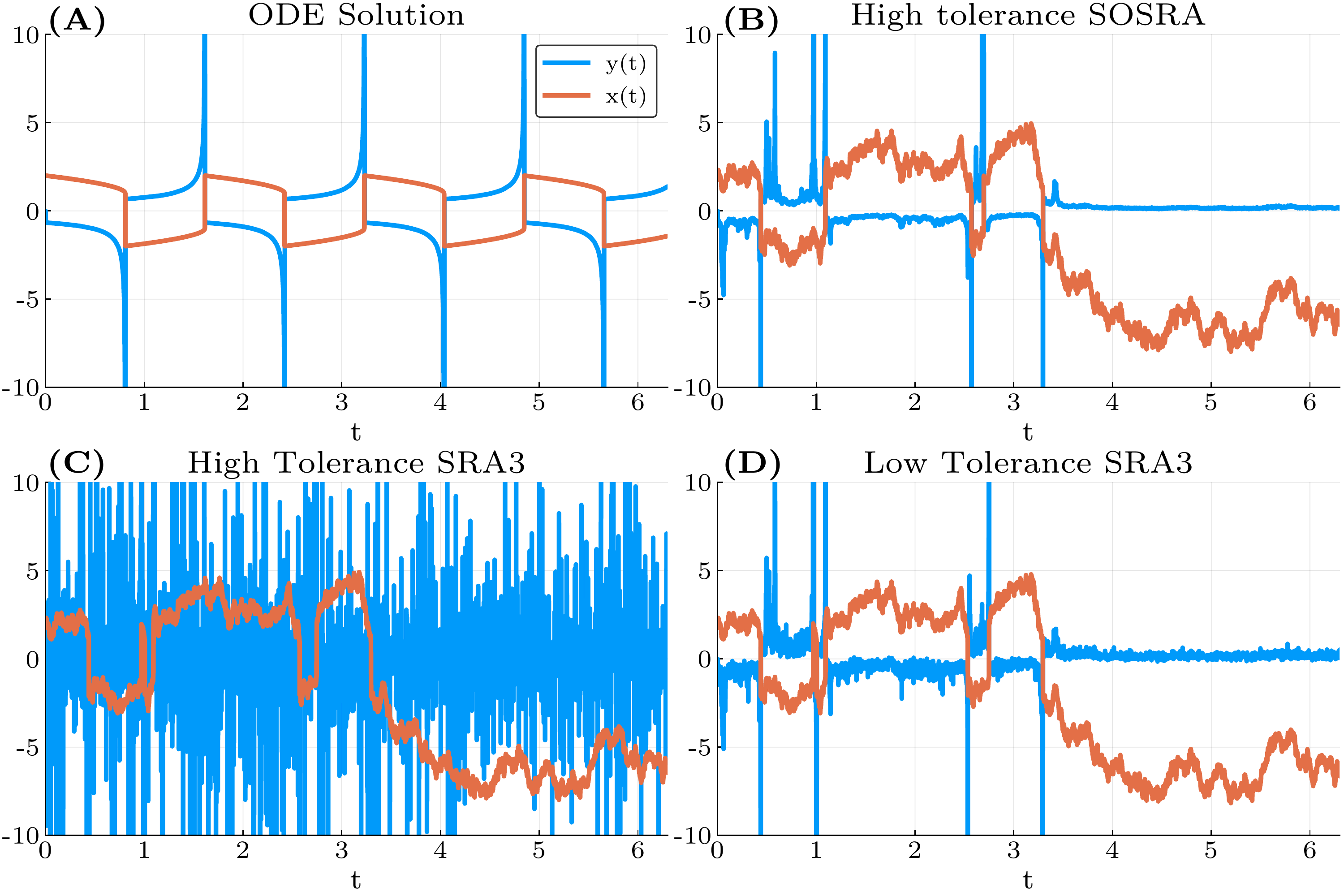}
			\par\end{centering}
		\caption{\textbf{Representative trajectories for solutions to the Van der Pol
				equations with additive noise.} Each of these trajectories are solved
			with the same underlying Brownian process. \textbf{(A)} The solution
			to the ODE with the explicit Runge-Kutta method Tsit5. \textbf{(B)
			}The solution to the SDE with tolerance $abstol=1,reltol=1/2^{1}$
			from SOSRA. \textbf{(C) }Solution to the SDE with tolerances $abstol=2^{-3},reltol=2^{-3}$
			with SRA3. \textbf{(D) }Solution to the SDE with tolerances $abstol=2^{-6},reltol=2^{-4}$
			with SRA3.\label{fig:AdditiveVanDerPol}}
	\end{figure}
	\par\end{center}

\subsubsection{Additive Van Der Pol Stiffness Detection}

In addition to testing efficiency, we used this to test the stiffness
detection in SOSRA2. Using a safety factor of $\omega=5$, we added
only two lines of code to make the algorithm print out the timings
for which the algorithm predicts stiffness. The results on two trajectories
were computed and are shown in Figure \ref{fig:Stiffness-detection-in-additive}.
The authors note that the stiffness detection algorithms are surprisingly
robust without any tweaking being done and are shown to not give almost
any false positives nor false negatives on this test problem. While
this safety factor is set somewhat high in comparison to traditional
ODE stiffness detection, we note that these algorithms were designed
to efficiently handle mild stiffness and thus we see it as a benefit
that they only declare stiffness when it appears to be in the regime
which is more suitable for implicit methods.
\begin{center}
	\begin{figure}
		\begin{centering}
			\includegraphics[scale=0.55]{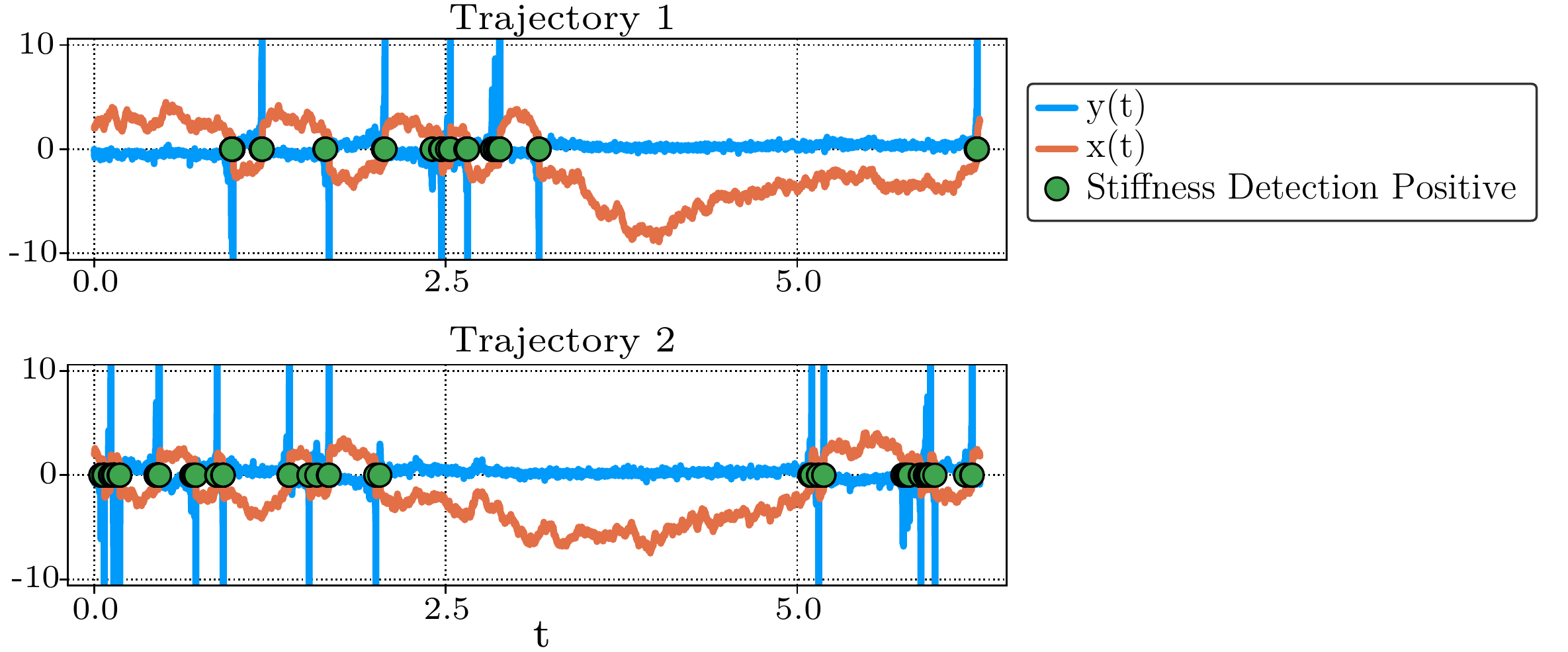}
			\par\end{centering}
		\caption{\textbf{Stiffness detection in the Van der Pol equations with additive
				noise Equation \ref{eq:van_der_pol_add}.} Two representative trajectories
			to Equation \ref{fig:AdditiveVanDerPol} are plotted. The green dots
			indicate time-points where the stiffness detection algorithm detected
			stiffness. \label{fig:Stiffness-detection-in-additive}}
	\end{figure}
	\par\end{center}

\subsection{SOSRI Numerical Efficiency Experiments}

\subsubsection{Multiplicative Noise Lotka-Volterra (2 Non-Stiff SDEs)}

To test the efficiency we plotted a work-precision diagram with SRIW1,
SOSRI, SOSRI2, and the fixed time step Euler-Maruyama and a Runge-Kutta
Milstein schemes for Equation \ref{eq:ex1} and the multiplicative
noise Lotka-Volterra Equation \ref{eq:lotka_mult}. As with Equation
\ref{eq:lotka_add}, Equation \ref{eq:lotka_mult} does not have an
analytical solution so a reference solution was computed using a low
tolerance solution via SOSRI for each Brownian trajectory. The results
show that there is a minimal difference in efficiency between the
SRI algorithms for errors over the interval $\left[10^{-6},10^{-2}\right]$,
while these algorithms are all significantly more efficient than the
lower order algorithms when the required error is $<10^{-2}$ (Figure
\ref{fig:SOSRI-efficiency}A-D). The weak error work-precision diagrams
show that when using between 100 to 10,000 trajectories, the weak
error is less than the sample error in the regime where there is no
discernible efficiency difference between the SRI methods.These results
show that in the regime of mild accuracy on non-stiff equations, these
methods are much more efficient than low order methods yet achieve
the same efficiency as the non-stability optimized SRI variants. Note
that these results also show the conservativeness of the error estimators.
\begin{center}
	\begin{figure}
		\begin{centering}
			\includegraphics[scale=0.33]{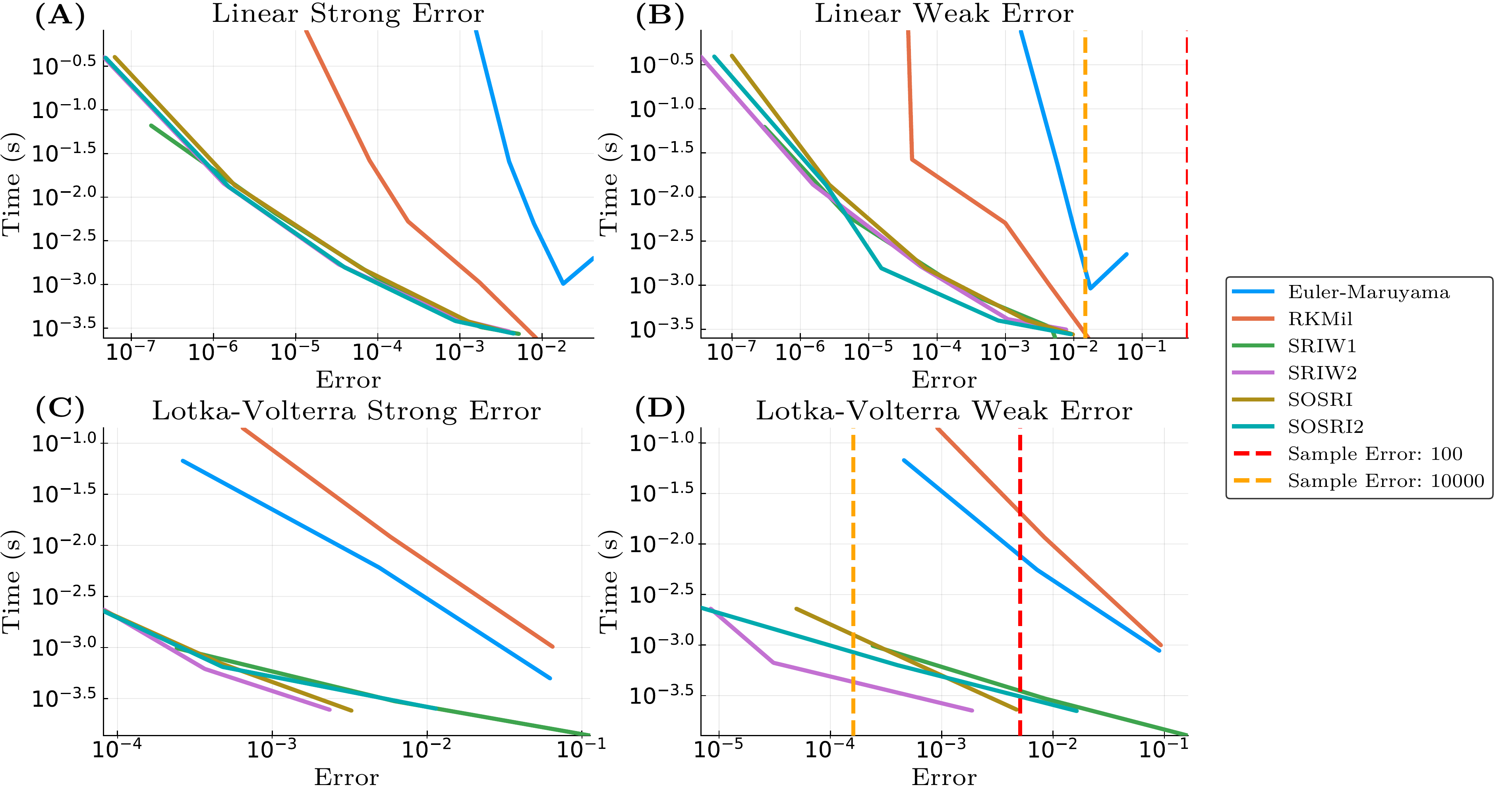}
			\par\end{centering}
		\caption{\textbf{SOSRI efficiency on non-stiff test equations.} The error was
			taken as the average of 10,000 trajectories for the Equation \ref{eq:ex1}
			and 100 trajectories for the Lokta-Volterra Equation \ref{eq:lotka_mult}.
			The sample error was determined for the weak error as the normal 95\%
			confidence interval for the mean using the variance of the true solution
			Equation \ref{eq:sol1} or the variance of the estimated true solutions
			via low tolerance solutions. The time is the average time to compute
			a trajectory and is averaged over 1000 runs at the same tolerance
			or step size.\textbf{ (A)} Shown are the work-precision plots for
			the methods on Equation \ref{eq:ex1}. Each of the adaptive time-stepping
			methods solved the problem on the interval using changing values of
			tolerances, with $tol=abstol=reltol$ starting at $10^{-1}$ and ending
			at $10^{-5}$ going in increments of $10$. The fixed time-stepping
			methods used time steps of size $h=5^{-2}$ to $h=5^{-7}$, changing
			the value by factors of $5$. The error is the strong $l_{2}$ error
			computed over the time series. \textbf{(B)} Same setup as the previous
			plot but using the weak error at the final time-point. \textbf{(C)
			}Shown are the work-precision plots for the methods on the multiplicative
			noise Lotka-Volterra Equation \ref{eq:lotka_mult}. Each of the adaptive
			time-stepping methods solved the problem on the interval using changing
			values of tolerances, with $tol=abstol=reltol$ starting at $4^{-2}$
			and ending at $4^{-4}$ going in increments of $4$. The fixed time-stepping
			methods used time steps of size $h=1/12^{-2.5}$ to $h=1/12^{-6.5}$,
			changing the value by factors of $12$. The error is the strong $l_{2}$
			error computed over the time series. \textbf{(D)} Same setup as the
			previous plot but using the weak error at the final time-point. \label{fig:SOSRI-efficiency}}
	\end{figure}
	\par\end{center}

\subsubsection{Epithelial-Mesenchymal Transition (EMT) Model (20 Pathwise Stiff
	SDEs) }

To test the real consequences of the enhanced stability, we use the
Epithelial-Mesenchymal Transition (EMT) model of 20 pathwise stiff
reaction equations introduced in \cite{RN3358}, studied as a numerical
test in \cite{RN3787}, and written in Section \ref{subsec:Epithelial-Mesenchymal-Transitio}.
In the previous work it was noted that $t\in\left[0,1\right]$ was
a less stiff version of this model. Thus we first tested the speed
that the methods could solve for 10,000 trajectories with no failures
due to numerical instabilities. The tolerances were tuned for each
method by factors of 2 and finding the largest values that were stable.
Since SOSRI demonstrated that its stability is much higher than even
SOSRI2, we show the effect of tolerance changes on SOSRI as well.
The results show that at similar tolerances the SOSRI method takes
nearly 5x less time than SRIW1 (Table \ref{tab:SRI-times-for-Oval2}).
However, there is an upper bound on the tolerances before the adaptivity
is no longer able to help keep the method stable. For SRIW1, this
bound is much lower, causing it to run more than 15x slower than the
fastest SOSRI setup. Interestingly SOSRI2 required a higher tolerance
than SRIW1 but was 3x faster than SRIW1's fastest setup. We note that
SOSRI's highest relative tolerance $2^{-7}\approx7\times10^{-3}$
is essentially requiring 4 digits of accuracy (in strong error) when
considering the conservativeness of the error estimator, which is
far beyond the accuracy necessary in many cases. Lastly, we note that
the SOSRI method is able to solve for 10,000 stable trajectories more
than 60x faster than any of the tested fixed time step methods.

\begin{center}
	\begin{table}
		\begin{centering}
			{\footnotesize
			\begin{tabular}{|c|c|c|c|c|}
				\hline
				Algorithm & Abstol & Reltol & Run-time (seconds) & Relative Time (vs SOSRI)\tabularnewline
				\hline
				\hline
				SOSRI & $2^{-7}$ & $2^{-4}$ & 2.62 & 1.0x\tabularnewline
				\hline
				SOSRI & $2^{-7}$ & $2^{-6}$ & 2.75 & 1.0x\tabularnewline
				\hline
				SOSRI & $2^{-12}$ & $2^{-15}$ & 8.78 & 3.3x\tabularnewline
				\hline
				SOSRI & $2^{-13}$ & $2^{-7}$ & 3.05 & 1.2x\tabularnewline
				\hline
				SOSRI2 & $2^{-12}$ & $2^{-15}$ & 8.69 & 3.3x\tabularnewline
				\hline
				SOSRI2 & $2^{-13}$ & $2^{-11}$ & 5.56 & 2.2x\tabularnewline
				\hline
				SRIW1 & $2^{-13}$ & $2^{-7}$ & 15.16 & 5.8x\tabularnewline
				\hline
				Euler-Maruyama &  &  & 169.96 & 64.8x\tabularnewline
				\hline
				Runge-Kutta Milstein &  &  & 182.59 & 69.6x\tabularnewline
				\hline
				Fixed Time-step SRIW1 &  &  & 424.30 & 161.7x\tabularnewline
				\hline
				DISTM $\left(\theta=\frac{1}{2}\right)$ &  &  & 8912.91 & 3396x\tabularnewline
				\hline
			\end{tabular}
			}%
			\par\end{centering}
		\caption{\textbf{SRI times for the the EMT model on $t\in\left[0,1\right]$.}
			The equations were solved 10,000 times with the given tolerances to
			completion and the elapsed time was recorded. The fixed time step
			methods had their $\Delta t$ determined as the largest $\Delta t$
			in increments of powers of 2 that produced no unstable trajectories,
			as shown in \cite{RN3787}. DISTM is the Drift-Implicit Stochastic
			$\theta$-Method\label{tab:SRI-times-for-Oval2} }
	\end{table}
	\par\end{center}

We then timed the run time to solve 10 trajectories in the $t\in\left[0,500\right]$
case (Table \ref{tab:SRI-times-for-Oval2-long}). This time we found
the optimal tolerance in terms of powers of $10$. Once again, SRIW1
needed a lower tolerance than is necessary in order to stay stable.
SOSRI is able to solve the problem only asking for around $tol=10^{-2}$,
while the others require more (especially in absolute tolerance as
there is a stiff reactant whose values travel close to zero). One
interesting point to note is that at similar tolerances both SOSRI
and SOSRI2 receive similar timings and both over 6 times faster than
the fastest SRIW1 tolerance setup. Both are nearly twice as fast as
SRIW1 when matching tolerances as well. Given the conservativeness
of the error estimators generally being around 2 orders of magnitude
more precise than the local error estimate, the low tolerance solutions
are accurate enough for many phenomenological experiments and thus
present a good speedup over previous methods. The timings for Euler-Maruyama
and Runge-Kutta Milstein schemes are omitted since the tests were
unable to finish. From the results of \cite{RN3787} we note that
the average $dt$ for SRIW1 on the edge of its stability had that
the smallest $dt$ was approximately $10^{-11}$. The stability region
for fixed step-size Euler-Maruyama is strictly smaller than SRIW1
(Figure \ref{fig:SOSRI-Stability-Regions.}) which suggests that it
would require around $5\times10^{12}$ time steps (with Runge-Kutta
Milstein being similar) to solve to $t=500$. Thus, given it takes
on our setup extrapolating the time given 170 seconds for $2^{20}$
steps, this projects to around $1.6\times10^{8}$ seconds, or approximately
5 years.
\begin{center}
	\begin{table}
		\begin{centering}
			{\footnotesize
			\begin{tabular}{|c|c|c|c|c|}
				\hline
				Algorithm & Abstol & Reltol & Run-time (seconds) & Relative Time (vs SOSRI)\tabularnewline
				\hline
				\hline
				SOSRI & $10^{-2}$ & $10^{-2}$ & 22.47 & 1.0x\tabularnewline
				\hline
				SOSRI & $10^{-4}$ & $10^{-4}$ & 73.62 & 3.3x\tabularnewline
				\hline
				SOSRI & $10^{-5}$ & $10^{-3}$ & 89.19 & 4.0x\tabularnewline
				\hline
				SOSRI2 & $10^{-4}$ & $10^{-4}$ & 76.12 & 3.4x\tabularnewline
				\hline
				SOSRI2 & $10^{-5}$ & $10^{-3}$ & 121.75 & 5.4x\tabularnewline
				\hline
				SRIW1 & $10^{-5}$ & $10^{-3}$ & 147.89 & 6.6x\tabularnewline
				\hline
				DIRKM $\left(\theta=\frac{1}{2}\right)$  &  &  & 7378.55 & 328.3x\tabularnewline
				\hline
				DIEM $\left(\theta=\frac{1}{2}\right)$ &  &  & 8796.47 & 391.4x\tabularnewline
				\hline
			\end{tabular}
			}%
			\par\end{centering}
		\caption{\textbf{SRI times for the the EMT model on $t\in\left[0,500\right]$.}
			The equations were solved 10 times with the given tolerances to completion
			and the elapsed time was recorded. The fixed timestep methods had
			their $\Delta t$ chosen by incrementing by $10^{-5}$ until 10 consecutive
			trajectories were stable. Drift-Implicit Euler Maruyama (DIEM) had
			$\Delta t=\frac{1}{60000}$ and Drift-Implicit Runge-Kutta Milstein (DIRKM)
			had $\Delta t=\frac{1}{50000}$.
			\label{tab:SRI-times-for-Oval2-long}}
	\end{table}
	\par\end{center}

\subsubsection{Retinoic Acid Stochastic Partial Differential Equation Model (6x20x100
	Semi-Stiff SDEs)}

As another test we applied the methods to a method of lines discretization
of a stochastic partial differential equation (SPDE) describing the
spatial regulation of the zebrafish hindbrain via retinoic acid signaling
( Section \ref{subsec:Retinoic-Acid-SPDE}) \cite{RN3805}. The discretization
results in a system of $6\times20\times100$ SDEs. Starting from an
initial zero state, a concentration gradient emerges over $t\in\left[0,500\right]$.
Each of the methods solved the problem at the highest tolerance that
was stable giving the results in Table \ref{tab:SRI-times-for-SPDE}.
Time stepping for this problem is heavily limited by the high diffusion
constant which results in a strict CFL condition for the 2nd order
finite difference discretization that is used (in the PDE sense),
making this problem's stepping stability-bound for explicit methods.
Because of this stiffness in the real axis, we found that the previous
high order adaptive method SRIW1 did not perform well on this problem
in comparison to Euler-Maruyama because the drift term is expensive
and the extra function calls outweighed the slightly larger timesteps.
However, the enhanced stability of the SOSRI and SOSRI2 methods allowed
for much larger time steps while keeping the same number of $f$ calls
per step, resulting in a more efficient solution when high accuracy
is not necessary. We note that the drift-implicit stochastic $\theta$-method
and drift implicit $\theta$Runge-Kutta Milstein methods were too
inefficient to estimate since their time steps were constrained to
be near that of the Euler-Maruyama equation due to divergence of the
Newton iterations. This SPDE could also be solved via SKenCarp by
using the transformation of Section \ref{sec:Affine-Transformation},
but from experiments on the PDE we note that efficient solution of
the implicit equations would require using a preconditioned Krylov
method due to the size of the system and thus it is left for future
investigation.
\begin{center}
	\begin{table}
		\begin{centering}
			{\footnotesize
			\begin{tabular}{|c|c|c|c|c|}
				\hline
				Algorithm & Abstol & Reltol & Run-time (seconds) & Relative Time (vs SOSRI)\tabularnewline
				\hline
				\hline
				SOSRI & $10^{-1}$ & $10^{-2}$ & 700.76 & 1.0x\tabularnewline
				\hline
				SOSRI2 & $10^{-3}$ & $10^{-3}$ & 1016.61 & 1.5x\tabularnewline
				\hline
				Euler-Maruyama &  &  & 1758.85 & 2.5x\tabularnewline
				\hline
				SRIW1 & $10^{-5}$ & $10^{-3}$ & 4205.52 & 6.0x\tabularnewline
				\hline
			\end{tabular}
			}%
			\par\end{centering}
		\caption{\textbf{SRI times for the the retinoic acid SPDE model on $t\in\left[0,500\right]$.}
			The equations were solved twice with the given tolerances to completion
			and the elapsed time was recorded. The tolerances were chosen as the
			highest pair of tolerances which did not diverge (going up by powers
			of 10). Note that none of the cases did the two timings vary by more
			than 1\% of the total run time. Euler-Maruyama used time steps of
			$\Delta t=1/20000$ since we note that at $\Delta t=1/10000$ approximately
			half of the trajectories (simulating 10) were unstable. \label{tab:SRI-times-for-SPDE}}
	\end{table}
	\par\end{center}

\section{Discussion}

In this work we derived stability-optimized SRK methods for additive
and diagonal noise equations, and used a transformation to allow the
additive noise methods to solve affine noise problems. Many other
equations can be reduced to the additive noise case as well using
the same means. Importantly, our derivation methods utilized heavy
computational tools in order to approximately optimize otherwise intractable
equations. This same method of derivation can easily be scaled up
to higher orders, and by incorporating the coefficients for higher
conditions, efficiency can be optimized as well by adding the norm
of the principle error coefficients to the optimization function.
The majority of the search was performed using global optimizers in
massive parallel using a hand-optimized CUDA kernel for the numerical
integral of the characteristic function, replacing man-hours with
core-hours and effectively optimizing the method. The clear next steps
are to find SRA and SRI methods with minimal error estimates and sensible
stability regions for the cases in which lower strong error matters,
and similar optimizations on SRK methods developed for small noise
problems. We note that high strong order methods were investigated
because of their better trajectory-wise convergence, allowing for
a more robust solution and error estimation since our application
to transiently pathwise stiff equations requires such properties.

In this work we also derived L-stable methods for additive (and thus
multiplicative and affine) noise equations, and computationally could
not find an A-B-L stable method. While our method does not prove that
no 2-stage A-B-L method exists, we have at least narrowed down its
possibility. Additionally an extension of a well-known ESDIRK method
to additive noise was developed. These ESDIRK methods have an extension
which allows for mass-matrices in the problem formulation. Using singular
mass matrices, these methods also present themselves as integrators
for a form of SDAEs with deterministic constraints. This method has
an implicit-explicit (IMEX) extension and the stochastic extension
was compatible with both tableaus. We showed that this IMEX version
of the method could numerically converge at order 2.0 on a test problem
(matching the other SRA methods), indicating that it may achieve the
sufficient condition. As an adaptive high order IMEX method, the ODE
version of the method is a common choice for large discretizations
of PDEs. Thus this method could present itself as a potentially automatic
and efficient option for discretizations of large affine noise SPDEs
by being able to use a low number of time steps while minimizing the
amount of work required to solve the implicit equation. We note that
adaptivity along with efficient stage predictors was required to be
more efficient than the common stochastic theta methods since divergence
of quasi-Newton steps can be common if care is not taken. After engineering
the method with all of the components together, the benchmark results
showed large efficiency gains over both the previous drift-implicit
and stability-optimized explicit methods. While previous literature
questioned the applicability of L-stable integrators to stochastic
differential equations due to high error in the slow variables \cite{RN3811},
our computations show that this analysis may be mislead by analyzing
strong order 0.5 methods. With our higher strong order methods we
see sufficiently accurate results on real stiff problems, and this
is greatly helped by time stepping adaptivity.

The main caveat for our methods is the restrictions on the form of
noise. While we have shown that an enlarged class of problems (affine
noise) can handled by the integrators for additive noise problems,
this is still a very special case in the scope of possible SDEs. Diagonal
noise is a much expanded scope but is still constrained, and our implicit
methods were only derived for the additive noise case. Further research
should focus on the expansion of this these techniques to high order
adaptive ESDIRK diagonal noise integrators. In addition, when $g$
is non-zero a ``diagonal noise'' problem over the complex plane
does not have diagonal noise (due to the mixing of real and complex
parts from complex multiplication, and reinterpretation as a $2n$
real system). Thus these methods are not applicable to problems defined
in the complex plane with complex Wiener processes. Development of
similar integrators for commutative noise problems could allow for
similar performance benefits on such problems and is a topic for future
research.

Additionally, we were not able to sufficiently improve the stability
along the noise axis with our explicit diagonal noise methods. However,
this is likely due to explicitness in the noise term. Recent research
has shown that step splitting which utilize a predicted step in the
diffusion calcuation can significantly improve the stability of a
method \cite{RN3815,RN3816}. Given this, we conjecture that a form
of predictor-correction, such as:

\begin{align}
X_{n+1}&=X_{n}+\sum_{i=1}^{s}\alpha_{i}f\left(t_{n}+c_{i}^{(0)}h,H_{i}^{(0)}\right)\\
&+\sum_{i=1}^{s}\left(\beta_{i}^{(1)}I_{(1)}+\beta_{i}^{(2)}\frac{I_{(1,1)}}{\sqrt{h}}+\beta_{i}^{(3)}\frac{I_{(1,0)}}{h}+\beta_{i}^{(4)}\frac{I_{(1,1,1)}}{h}\right)g\left(t_{n}+c_{i}^{(1)}h\right)\label{eq:update-1}
\end{align}
with stages
\begin{align}
\overline{H_{i}^{(1)}} & =X_{n}+\sum_{j=1}^{s}A_{ij}^{(1)}f\left(t_{n}+c_{j}^{(0)}h,H_{j}^{(0)}\right)h+\sum_{j=1}^{s}B_{ij}^{(1)}g\left(t_{n}+c_{j}^{(1)}h,H_{i}^{(1)}\right)\sqrt{h}\label{eq:stages-1}\\
H_{i}^{(0)} & =X_{n}+\sum_{j=1}^{s}A_{ij}^{(0)}f\left(t_{n}+c_{j}^{(0)}h,H_{j}^{(0)}\right)h+\sum_{j=1}^{s}B_{ij}^{(0)}g\left(t_{n}+c_{j}^{(1)}h,\overline{H_{i}^{(1)}}\right)\frac{I_{(1,0)}}{h}\\
H_{i}^{(1)} & =X_{n}+\sum_{j=1}^{s}A_{ij}^{(1)}f\left(t_{n}+c_{j}^{(0)}h,H_{j}^{(0)}\right)h+\sum_{j=1}^{s}B_{ij}^{(1)}g\left(t_{n}+c_{j}^{(1)}h,\overline{H_{i}^{(1)}}\right)\sqrt{h}\nonumber
\end{align}
could improve the noise stability of the method while keeping explicitness
and the same tableau. However, proper convergence and stability analysis
would require significant effort.

Our timings show that the current high order SRK methods are stability-bound
and that when scientific studies are only looking for small amounts
of accuracy in stochastic simulations, most of the computational effort
is lost to generating more accurate than necessary solutions in order
to satisfy stability constraints. For additive noise problems we were
able to obtain solutions about 5x-30x faster and for diagonal noise
approximately 6x than the current adaptive methods (SRA1, SRA3, SRIW1),
while common methods like Euler-Maruyama and Drift-Implicit $\theta$
Runge-Kutta Milstein were in many cases hundreds of times slower or
in many cases could not even finish. We have also shown that these
methods are very robust even at high tolerances and have a tendency
to produce the correct qualitative results on semi-stiff equations
(via plots) even when the user chosen accuracy is low. Given that
the required user input is minimal and work over a large range of
stiffness, we see these as very strong candidates for default general
purpose solvers for problem-solving environments such as MATLAB and
Julia since they can easily and efficiently produce results which
are sufficiently correct. Due to a choice in the optimization, the
SOSRA and SOSRA2 methods are not as efficient at low tolerances as
SRA3, so SRA3 should be used when high accuracy is necessary (on additive
or affine noise problems). However, in many cases like integrating
to find steady distributions of bistable parameter regimes or generating
trajectories of phonomenological models, this ability to quickly get
a more course estimate is valuable.

The stiffness detection in SDEs is a novel addition which we have
demonstrated can act very robustly. It has a control parameter $\omega$
which can be used to control the false positive and false negative
rate as needed. Note that stiff methods can achieve similar largest
eigenvalue estimates directly from the Jacobians of $f$ (and $g$)
given that the methods are implicit (or in the case of Rosenbrock
methods, the Jacobian must still be computed), and thus this can be
paired with a stiff solver to allow for automatic switching between
stiff and non-stiff solvers. Given that the cost for such stiffness
checks is minimal and the demonstrated efficiency of the implicit
methods on stiff equations, we are interested in future studies on
the efficiency of such composite method due to the stochastic nature
of stiffness in SDEs.

\appendix
\section{Appendix I: Test Equations}

\subsection{Additive Noise Test Equation}

\begin{equation}
dX_{t}=\left(\frac{\beta}{\sqrt{1+t}}-\frac{1}{2\left(1+t\right)}X_{t}\right)dt+\frac{\alpha\beta}{\sqrt{1+t}}dW_{t},\thinspace\thinspace\thinspace X_{0}=\frac{1}{2},\label{eq:ex3}
\end{equation}
where $\alpha=\frac{1}{10}$ and $\beta=\frac{1}{20}$ with true
solution
\begin{equation}
X_{t}=\frac{1}{\sqrt{1+t}}X_{0}+\frac{\beta}{\sqrt{1+t}}\left(t+\alpha W_{t}\right).\label{eq:sol3}
\end{equation}

\subsection{Diagonal Noise Test Equation}

\begin{equation}
dX_{t}=\alpha X_{t}dt+\beta X_{t}dW_{t}\thinspace\thinspace\thinspace X_{0}=\frac{1}{2},\label{eq:ex1}
\end{equation}
where $\alpha=\frac{1}{10}$ and $\beta=\frac{1}{20}$ with true solution
\begin{equation}
X_{t}=X_{0}e^{\left(\beta-\frac{\alpha^{2}}{2}\right)t+\alpha W_{t}}.\label{eq:sol1}
\end{equation}

\subsection{Split Additive Test Equation}

\begin{equation}
dX_{t}=\left(f_{1}(t,X_{t})+f_{2}(t,X_{t})\right)dt+\frac{\alpha\beta}{\sqrt{1+t}}dW_{t},\thinspace\thinspace\thinspace X_{0}=\frac{1}{2},\label{eq:ex3-1}
\end{equation}
where
\begin{align*}
f_{1}(t,X_{t}) & =\frac{\beta}{\sqrt{1+t}}\\
f_{2}(t,X_{t}) & =-\frac{1}{2\left(1+t\right)}X_{t}
\end{align*}
with true solution Equation \ref{eq:sol3}.

\subsection{Additive Noise Lotka-Volterra}

\begin{align}
dx & =\left(ax-bxy\right)dt+\sigma_{A}dW_{t}^{1}\nonumber \\
dy & =\left(-cy+dxy\right)dt+\sigma_{A}dW_{t}^{2}\label{eq:lotka_add}
\end{align}
where $a=1.5$, $b=1$, $c=3.0$, $d=1.0$, $\sigma_{A}=0.01$.

\subsection{Additive Noise Van Der Pol}

The driven Van Der Pol equation is
\begin{align}
dy & =\mu((1-x^{2})y-x)dt\nonumber \\
dx & =ydt\label{eq:van_der_pol}
\end{align}
The additive noise variant is
\begin{align}
dy & =\mu((1-x^{2})y-x)dt+\rho dW_{t}^{(1)}\nonumber \\
dx & =y+\rho dW_{t}^{(2)}\label{eq:van_der_pol_add}
\end{align}

\subsection{Multiplicative Noise Lotka-Volterra}

\begin{align}
dx & =\left(ax-bxy\right)dt+\sigma_{A}dW_{t}^{1}\nonumber \\
dy & =\left(-cy+dxy\right)dt+\sigma_{A}dW_{t}^{2}\label{eq:lotka_mult}
\end{align}
where $a=1.5$, $b=1$, $c=3.0$, $d=1.0$, $\sigma_{A}=0.01$.

\subsection{Epithelial-Mesenchymal Transition Model \label{subsec:Epithelial-Mesenchymal-Transitio}}

The Epithelial-\\
Mesenchymal Transition (EMT) model is given by the
following system of SDEs which correspond to a chemical reaction network
modeled via mass-action kinetics with Hill functions for the feedbacks.
This model was introduced in \cite{RN3358}.

{\footnotesize{}\allowdisplaybreaks
	\begin{eqnarray*}
		A & = & \left(\left(\left[TGF\right]+\left[TGF0\right](t)\right)/J0_{snail}\right)^{n0_{snail}}+\left(\left[OVOL2\right]/J1_{snail}\right)^{n1_{snail}}\\
		\frac{d\left[snail1\right]_{t}}{dt} & = & k0_{snail}+k_{snail}\frac{\left(\left(\left[TGF\right]+\left[TGF0\right](t)\right)/J0_{snail}\right)^{n0_{snail}}}{\left(1+A\right)\left(1+\left[SNAIL\right]/J2_{snail}\right)}\\
		& - & kd_{snail}\left(\left[snail1\right]-\left[SR\right]\right)-kd_{SR}\left[SR\right]\\
		\frac{d\left[SNAIL\right]}{dt} & = & k_{SNAIL}\left(\left[snail1\right]-\left[SR\right]\right)-kd_{SNAIL}\left[SNAIL\right]\\
		\frac{d\left[miR34\right]}{dt} & = & kO_{34}+\frac{k_{34}}{1+\left(\left[SNAIL\right]/J1_{34}\right)^{n1_{34}}+\left(\left[ZEB\right]/J2_{34}\right)^{n2_{34}}}\\
		& - & kd_{34}\left(\left[miR34\right]-\left[SR\right]\right)-\left(1-\lambda_{SR}\right)kd_{SR}\left[SR\right]\\
		\frac{d\left[SR\right]}{dt} & = & Tk\left(K_{SR}\left(\left[snail1\right]-\left[SR\right]\right)\left(\left[miR34\right]-\left[SR\right]\right)-\left[SR\right]\right)\\
		\frac{d\left[zeb\right]}{dt} & = & k0_{zeb}+k_{zeb}\frac{\left(\left[SNAIL\right]/J1_{zeb}\right)^{n1_{zeb}}}{1+\left(\left[SNAIL\right]/J1_{zeb}\right)^{n1_{zeb}}+\left(\left[OVOL2\right]/J2_{zeb}\right)^{n2_{zeb}}}\\
		& - & kd_{zeb}\left(\left[zeb\right]-\sum_{i=1}^{5}C_{5}^{i}\left[ZR\right]\right)-\sum_{i=1}^{5}kd_{ZR_{i}}C_{5}^{i}\left[ZR_{i}\right]\\
		\frac{d\left[ZEB\right]}{dt} & = & k_{ZEB}\left(\left[zeb\right]-\sum_{i=1}^{5}C_{5}^{i}\left[ZR_{i}\right]\right)-kd_{ZEB}\left[ZEB\right]\\
		\frac{d\left[miR200\right]}{dt} & = & k0_{200}+\frac{k_{200}}{1+\left(\left[SNAIL\right]/J1_{200}\right)^{n1_{200}}+\left(\left[ZEB\right]/J2_{200}\right)^{n2_{200}}}\\
		& - & kd_{200}\left(\left[miR200\right]-\sum_{i=1}^{5}iC_{5}^{i}\left[ZR_{i}\right]-\left[TR\right]\right)\\
		& - & \sum_{i=1}^{5}\left(1-\lambda_{i}\right)kd_{ZR_{i}}C_{5}^{i}i\left[ZR_{i}\right]-\left(1-\lambda_{TR}\right)kd_{TR}\left[TR\right]\\
		\frac{d\left[ZR_{1}\right]}{dt} & = & Tk\left(K_{1}\left(\left[miR200\right]-\sum_{i=1}^{5}iC_{5}^{i}\left[ZR_{i}\right]-\left[TR\right]\right)\right.\\
		&  & \left.\left(\left[zeb\right]-\sum_{i=1}^{5}C_{5}^{i}\left[ZR_{i}\right]\right)-\left[ZR_{1}\right]\right)\\
		\frac{d\left[ZR_{2}\right]}{dt} & = & Tk\left(K_{2}\left(\left[miR200\right]-\sum_{i=1}^{5}iC_{5}^{i}\left[ZR_{i}\right]-\left[TR\right]\right)\left[ZR_{1}\right]-\left[ZR_{2}\right]\right)\\
		\frac{d\left[ZR_{3}\right]}{dt} & = & Tk\left(K_{3}\left(\left[miR200\right]-\sum_{i=1}^{5}iC_{5}^{i}\left[ZR_{i}\right]-\left[TR\right]\right)\left[ZR_{1}\right]-\left[ZR_{3}\right]\right)\\
		\frac{d\left[ZR_{4}\right]}{dt} & = & Tk\left(K_{4}\left(\left[miR200\right]-\sum_{i=1}^{5}iC_{5}^{i}\left[ZR_{i}\right]-\left[TR\right]\right)\left[ZR_{1}\right]-\left[ZR_{4}\right]\right)\\
		\frac{d\left[ZR_{5}\right]}{dt} & = & Tk\left(K_{5}\left(\left[miR200\right]-\sum_{i=1}^{5}iC_{5}^{i}\left[ZR_{i}\right]-\left[TR\right]\right)\left[ZR_{1}\right]-\left[ZR_{5}\right]\right)\\
		\frac{d\left[tgf\right]}{dt} & = & k_{tgf}-kd_{tgf}\left(\left[tgf\right]-\left[TR\right]\right)-kd_{TR}\left[TR\right]\\
		\frac{d\left[TGF\right]}{dt} & = & k0_{TGF}+k_{TGF}\left(\left[tgf\right]-\left[TR\right]\right)-kd_{TGF}\left[TGF\right]\\
		\frac{d\left[TR\right]}{dt} & = & Tk\left(K_{TR}\left(\left[miR200\right]-\sum_{i=1}^{5}iC_{5}^{i}\left[ZR_{i}\right]-\left[TR\right]\right)\left(\left[tgf\right]-\left[TR\right]\right)-\left[TR\right]\right)\\
		\frac{d\left[Ecad\right]}{dt} & = & k0_{E}+\frac{k_{E1}}{1+\left(\left[SNAIL\right]/J1_{E}\right)^{n1_{E}}}+\frac{k_{E2}}{1+\left(\left[ZEB\right]/J2_{E}\right)^{n2_{E}}}-kd_{E}\left[Ecad\right]\\
		B & = & k_{V1}\frac{\left(\left[SNAIL\right]/J1_{V}\right)^{n1_{V}}}{1+\left(\left[SNAIL\right]/J1_{V}\right)^{n1_{V}}}+k_{V2}\frac{\left(\left[ZEB\right]/J2_{V}\right)^{n2_{V}}}{1+\left(\left[ZEB\right]/J2_{V}\right)^{n2_{V}}}\\
		\frac{d\left[Vim\right]}{dt} & = & k0_{V}+\frac{B}{\left(1+\left[OVOL2\right]/J3_{V}\right)}-kd_{V}\left[Vim\right]\\
		\frac{d\left[OVOL2\right]}{dt} & = & k0_{0}+k_{0}\frac{1}{1+\left(\left[ZEB\right]/J_{0}\right)^{n_{0}}}-kd_{O}\left[OVOL2\right]\\
		\frac{d\left[OVOL2\right]_{p}}{dt} & = & k_{Op}\left[OVOL2\right]-kd_{Op}\left[OVOL2\right]_{p}
	\end{eqnarray*}
} where{\footnotesize{}
\begin{eqnarray*}
	\sum_{i=1}^{5}iC_{5}^{i}\left[ZR_{i}\right] & = & 5\left[ZR_{1}\right]+20\left[ZR_{2}\right]++30\left[ZR_{3}\right]+20\left[ZR_{4}\right]+5\left[ZR_{5}\right],\\
	\sum_{i=1}^{5}C_{5}^{i}\left[ZR_{i}\right] & = & 5\left[ZR_{1}\right]+10\left[ZR_{2}\right]+10\left[ZR_{3}\right]+5\left[ZR_{4}\right]+\left[ZR_{5}\right],\\
	\left[TGF0\right](t) & = & \begin{cases}
		\frac{1}{2} & t>100\\
		0 & o.w.
	\end{cases}
\end{eqnarray*}
}{\footnotesize \par}

The parameter values are given in Table \ref{tab:params}.

\begin{table}
	\centering{}{\footnotesize{} }%
	\begin{tabular}{|c|c|c|c|c|c|c|c|}
		\hline
		{\footnotesize{}Parameter } & {\footnotesize{}Value } & {\footnotesize{}Parameter } & {\footnotesize{}Value } & {\footnotesize{}Parameter } & {\footnotesize{}Value } & {\footnotesize{}Parameter } & {\footnotesize{}Value}\tabularnewline
		\hline
		\hline
		{\footnotesize{}$J1_{200}$ } & {\footnotesize{}3 } & {\footnotesize{}$J1_{E}$ } & {\footnotesize{}0.1 } & {\footnotesize{}$K_{2}$ } & {\footnotesize{}1 } & {\footnotesize{}$k0_{O}$ } & {\footnotesize{}0.35}\tabularnewline
		\hline
		{\footnotesize{}$J2_{200}$ } & {\footnotesize{}0.2 } & {\footnotesize{}$J2_{E}$ } & {\footnotesize{}0.3 } & {\footnotesize{}$K_{3}$ } & {\footnotesize{}1 } & {\footnotesize{}$kO_{200}$ } & {\footnotesize{}0.0002}\tabularnewline
		\hline
		{\footnotesize{}$J1_{34}$ } & {\footnotesize{}0.15 } & {\footnotesize{}$J1_{V}$ } & {\footnotesize{}0.4 } & {\footnotesize{}$K_{4}$ } & {\footnotesize{}1 } & {\footnotesize{}$kO_{34}$ } & {\footnotesize{}0.001}\tabularnewline
		\hline
		{\footnotesize{}$J2_{34}$ } & {\footnotesize{}0.35 } & {\footnotesize{}$J2_{V}$ } & {\footnotesize{}0.4 } & {\footnotesize{}$K_{5}$ } & {\footnotesize{}1 } & {\footnotesize{}$kd_{snail}$ } & {\footnotesize{}0.09}\tabularnewline
		\hline
		{\footnotesize{}$J_{O}$ } & {\footnotesize{}0.9 } & {\footnotesize{}$J3_{V}$ } & {\footnotesize{}2 } & {\footnotesize{}$K_{TR}$ } & {\footnotesize{}20 } & {\footnotesize{}$kd_{tgf}$ } & {\footnotesize{}0.1}\tabularnewline
		\hline
		{\footnotesize{}$J0_{snail}$ } & {\footnotesize{}0.6 } & {\footnotesize{}$J1_{zeb}$ } & {\footnotesize{}3.5 } & {\footnotesize{}$K_{SR}$ } & {\footnotesize{}100 } & {\footnotesize{}$kd_{zeb}$ } & {\footnotesize{}0.1}\tabularnewline
		\hline
		{\footnotesize{}$J1_{snail}$ } & {\footnotesize{}0.5 } & {\footnotesize{}$J2_{zeb}$ } & {\footnotesize{}0.9 } & {\footnotesize{}$TGF0$ } & {\footnotesize{}0 } & {\footnotesize{}$kd_{TGF}$ } & {\footnotesize{}0.9}\tabularnewline
		\hline
		{\footnotesize{}$J2_{snail}$ } & {\footnotesize{}1.8 } & {\footnotesize{}$K_{1}$ } & {\footnotesize{}1 } & {\footnotesize{}$Tk$ } & {\footnotesize{}1000 } & {\footnotesize{}$kd_{ZEB}$ } & {\footnotesize{}1.66}\tabularnewline
		\hline
		{\footnotesize{}$k0_{snail}$ } & {\footnotesize{}0.0005 } & {\footnotesize{}$k0_{zeb}$ } & {\footnotesize{}0.003 } & {\footnotesize{}$\lambda_{1}$ } & {\footnotesize{}0.5 } & {\footnotesize{}$k0_{TGF}$ } & {\footnotesize{}1.1}\tabularnewline
		\hline
		{\footnotesize{}$n1_{200}$ } & {\footnotesize{}3 } & {\footnotesize{}$n1_{snail}$ } & {\footnotesize{}2 } & {\footnotesize{}$\lambda_{2}$ } & {\footnotesize{}0.5 } & {\footnotesize{}$k0_{E}$ } & {\footnotesize{}5}\tabularnewline
		\hline
		{\footnotesize{}$n2_{200}$ } & {\footnotesize{}2 } & {\footnotesize{}$n1_{E}$ } & {\footnotesize{}2 } & {\footnotesize{}$\lambda_{3}$ } & {\footnotesize{}0.5 } & {\footnotesize{}$k0_{V}$ } & {\footnotesize{}5}\tabularnewline
		\hline
		{\footnotesize{}$n1_{34}$ } & {\footnotesize{}2 } & {\footnotesize{}$n2_{E}$ } & {\footnotesize{}2 } & {\footnotesize{}$\lambda_{4}$ } & {\footnotesize{}0.5 } & {\footnotesize{}$k_{E1}$ } & {\footnotesize{}15}\tabularnewline
		\hline
		{\footnotesize{}$n2_{34}$ } & {\footnotesize{}2 } & {\footnotesize{}$n1_{V}$ } & {\footnotesize{}2 } & {\footnotesize{}$\lambda_{5}$ } & {\footnotesize{}0.5 } & {\footnotesize{}$k_{E2}$ } & {\footnotesize{}5}\tabularnewline
		\hline
		{\footnotesize{}$n_{O}$ } & {\footnotesize{}2 } & {\footnotesize{}$n2_{V}$ } & {\footnotesize{}2 } & {\footnotesize{}$\lambda_{SR}$ } & {\footnotesize{}0.5 } & {\footnotesize{}$k_{V1}$ } & {\footnotesize{}2}\tabularnewline
		\hline
		{\footnotesize{}$n0_{snail}$ } & {\footnotesize{}2 } & {\footnotesize{}$n2_{zeb}$ } & {\footnotesize{}6 } & {\footnotesize{}$\lambda_{TR}$ } & {\footnotesize{}0.5 } & {\footnotesize{}$k_{V2}$ } & {\footnotesize{}5}\tabularnewline
		\hline
		{\footnotesize{}$k_{O}$ } & {\footnotesize{}1.2 } & {\footnotesize{}$k_{200}$ } & {\footnotesize{}0.02 } & {\footnotesize{}$k_{34}$ } & {\footnotesize{}0.01 } & {\footnotesize{}$k_{tgf}$ } & {\footnotesize{}0.05}\tabularnewline
		\hline
		{\footnotesize{}$k_{zeb}$ } & {\footnotesize{}0.06 } & {\footnotesize{}$k_{TGF}$ } & {\footnotesize{}1.5 } & {\footnotesize{}$k_{SNAIL}$ } & {\footnotesize{}16 } & {\footnotesize{}$k_{ZEB}$ } & {\footnotesize{}16}\tabularnewline
		\hline
		{\footnotesize{}$kd_{ZR_{1}}$ } & {\footnotesize{}0.5 } & {\footnotesize{}$kd_{ZR_{2}}$ } & {\footnotesize{}0.5 } & {\footnotesize{}$kd_{ZR_{3}}$ } & {\footnotesize{}0.5 } & {\footnotesize{}$kd_{ZR_{4}}$ } & {\footnotesize{}0.5}\tabularnewline
		\hline
		{\footnotesize{}$kd_{ZR_{5}}$ } & {\footnotesize{}0.5 } & {\footnotesize{}$kd_{O}$ } & {\footnotesize{}1.0 } & {\footnotesize{}$kd_{200}$ } & {\footnotesize{}0.035 } & {\footnotesize{}$kd_{34}$ } & {\footnotesize{}0.035}\tabularnewline
		\hline
		{\footnotesize{}$kd_{SR}$ } & {\footnotesize{}0.9 } & {\footnotesize{}$kd_{E}$ } & {\footnotesize{}0.05 } & {\footnotesize{}$kd_{V}$ } & {\footnotesize{}0.05 } & $k_{Op}$ & 10\tabularnewline
		\hline
		&  &  &  &  &  & $kd_{Op}$ & 10\tabularnewline
		\hline
	\end{tabular} \caption{Table of Parameter Values for the EMT Model.\label{tab:params}}
\end{table}

\subsection{Retinoic Acid SPDE Model \label{subsec:Retinoic-Acid-SPDE}}

{\tiny{}
	\begin{align*}
	d\left[RA_{out}\right] & =\left(\beta(x)+D\Delta\left[RA_{out}\right]-b\left[RA_{out}\right]+c\left[RA_{in}\right]\right)dt+\sigma_{RA_{out}}dW_{t}^{out}\\
	d\left[RA_{in}\right] & =\left(b\left[RA_{out}\right]+\delta\left[BP\right]\left[RA-RAR\right]-\left(\gamma\left[BP\right]+\eta+\frac{\alpha\left[RA-RAR\right]}{\omega+\left[RA-RAR\right]}-c\right)\left[RA_{in}\right]\right)dt\\
	d\left[RA-BP\right] & =\left(\gamma\left[BP\right]\left[RA_{in}\right]+\lambda\left[BP\right]\left[RA-RAR\right]-\left(\delta+\nu\left[RAR\right]\right)\left[RA-BP\right]\right)dt\\
	d\left[RA-RAR\right] & =\left(\nu\left[RA-BP\right]\left[RAR\right]-\lambda\left[BP\right]\left[RA-RAR\right]\right)dt+\sigma_{RA-RAR}\left[RA-RAR\right]dW_{t}^{RA-RAR}\\
	d\left[BP\right] & =\left(a-\lambda\left[BP\right]\left[RA-RAR\right]-\gamma\left[BP\right]\left[RA_{in}\right]+\left(\delta+\nu\left[RAR\right]\right)\left[RA-BP\right]-u\left[BP\right]+\frac{d\left[RA-RAR\right]}{e+\left[RA-RAR\right]}\right)dt\\
	d\left[RAR\right] & =\left(\zeta-\nu\left[RA-BP\right]\left[RAR\right]+\lambda\left[BP\right]\left[RA-RAR\right]-r\left[RAR\right]\right)dt
	\end{align*}
}{\scriptsize \par}

where $\beta(x)=\beta_{0}H(x-40)$ with $H$ the Heaviside step function
and $x=40$ is the edge of retinoic acid production \cite{RN3805}.
The space was chosen as $\left[-100,400\right]\times\left[0,100\right]$
with $\Delta x=\Delta y=5$. The boundary conditions were no-flex
on every side except the right side which had leaky boundary conditions
with parameter $kA=0.002$, though full no-flux does not noticably
change the results. The parameter values are given in Table \ref{tab:params2}.

\begin{table}
	\centering{}{\footnotesize{} }%
	\begin{tabular}{|c|c|c|c|c|c|}
		\hline
		{\footnotesize{}Parameter } & {\footnotesize{}Value } & {\footnotesize{}Parameter } & {\footnotesize{}Value } & {\footnotesize{}Parameter } & {\footnotesize{}Value }\tabularnewline
		\hline
		\hline
		$\sigma_{RA_{in}}$,$\sigma_{RA-RAR}$,$\sigma_{RA_{out}}$ & 0.1 & $\omega$ & 100 & $u$ & 0.01\tabularnewline
		\hline
		b & 0.17 & $\gamma$ & 3.0 & $d$ & 0.1\tabularnewline
		\hline
		$\alpha$ & 10000  & $\delta$ & 0.0013  & $e$ & 1\tabularnewline
		\hline
		$\beta_{0}$ & 1 & $\eta$ & 0.0001  & $a$ & 1\tabularnewline
		\hline
		$c$ & 0.1 & $r$ & 0.0001  & $\zeta$ & 0.02\tabularnewline
		\hline
		$\nu$ & 0.85 & $\lambda$ & 0.85 & $D$ & 250.46\tabularnewline
		\hline
	\end{tabular} \caption{Table of Parameter Values for the EMT Model.\label{tab:params2}}
\end{table}

\section{Appendix I: SKenCarp, SOSRA, and SOSRI Tableaus}

All entries not listed are zero.

\subsection{SKenCarp Exact Values\label{subsec:SKenCarpExact}}

{\small
\begin{align*}
	K_1 & = 87294609440832483406992237 \\
	K_2 & = -53983406399371387722712393713535786276 \\
	K_3 & = 26826820\sqrt{6853072660943221216270384658311461343029149665543510113394397} \\
	K_4 & =\frac{K_1\left(K_2-K_3\right)}{4868738516734691891458097}\\
	B_{2,1}^{(0)} & =\frac{K_4-354038415192410790619483213666362001932}{210758174113231167877981435258781706648},\\
	B_{4,3}^{(0)} & =\frac{K_2-K_3}{8606625878152317177894269252900546591},\\
	B_{i,j}^{(0)} & =0\thinspace o.w.
\end{align*}
}%

\subsection{SOSRA \label{subsec:SOSRA}}
\begin{center}
	\begin{tabular}{|c|c||c|c|}
		\hline
		Coefficient & Value & Coefficient & Value\tabularnewline
		\hline
		\hline
		$\alpha_{1}$ & 0.2889874966892885 & $\beta_{3}^{(1)}$ & 0.27753845684143835\tabularnewline
		\hline
		$\alpha_{2}$ & 0.6859880440839937 & $\beta_{1}^{(2)}$ & 0.4237535769069274\tabularnewline
		\hline
		$\alpha_{3}$ & 0.025024459226717772 & $\beta_{2}^{(2)}$ & 0.6010381474428539\tabularnewline
		\hline
		$c_{1}^{(0)}$ & 0 & $\beta_{3}^{(2)}$ & -1.0247917243497813\tabularnewline
		\hline
		$c_{2}^{(0)}$ & 0.6923962376159507 & $A_{2,1}^{(0)}$ & 0.6923962376159507\tabularnewline
		\hline
		$c_{3}^{(0)}$ & 1 & $A_{3,1}^{(0)}$ & -3.1609142252828395\tabularnewline
		\hline
		$c_{1}^{(1)}$ & 0 & $A_{3,2}^{(0)}$ & 4.1609142252828395\tabularnewline
		\hline
		$c_{2}^{(1)}$ & 0.041248171110700504 & $B_{2,1}^{(0)}$ & 1.3371632704399763\tabularnewline
		\hline
		$c_{3}^{(1)}$ & 1 & $B_{3,1}^{(0)}$ & 1.442371048468624\tabularnewline
		\hline
		$\beta_{1}^{(1)}$ & -16.792534242221663 & $B_{3,2}^{(0)}$ & 1.8632741501139225\tabularnewline
		\hline
		$\beta_{2}^{(1)}$ & 17.514995785380226 &  & \tabularnewline
		\hline
	\end{tabular}
	\par\end{center}

\subsection{SOSRA2 \label{subsec:SOSRA2}}
\begin{center}
	\begin{tabular}{|c|c||c|c|}
		\hline
		Coefficient & Value & Coefficient & Value\tabularnewline
		\hline
		\hline
		$\alpha_{1}$ & 0.4999999999999998 & $\beta_{3}^{(1)}$ & 0.07561967854316998\tabularnewline
		\hline
		$\alpha_{2}$ & -0.9683897375354181 & $\beta_{1}^{(2)}$ & 1\tabularnewline
		\hline
		$\alpha_{3}$ & 1.4683897375354185 & $\beta_{2}^{(2)}$ & -0.8169981105823436\tabularnewline
		\hline
		$c_{1}^{(0)}$ & 0 & $\beta_{3}^{(2)}$ & -0.18300188941765633\tabularnewline
		\hline
		$c_{2}^{(0)}$ & 1 & $A_{2,1}^{(0)}$ & 1\tabularnewline
		\hline
		$c_{3}^{(0)}$ & 1 & $A_{3,1}^{(0)}$ & 0.9511849235504364\tabularnewline
		\hline
		$c_{1}^{(1)}$ & 0 & $A_{3,2}^{(0)}$ & 0.04881507644956362\tabularnewline
		\hline
		$c_{2}^{(1)}$ & 1 & $B_{2,1}^{(0)}$ & 0.7686101171003622\tabularnewline
		\hline
		$c_{3}^{(1)}$ & 1 & $B_{3,1}^{(0)}$ & 0.43886792994934987\tabularnewline
		\hline
		$\beta_{1}^{(1)}$ & 0 & $B_{3,2}^{(0)}$ & 0.7490415909204886\tabularnewline
		\hline
		$\beta_{2}^{(1)}$ & 0.92438032145683 &  & \tabularnewline
		\hline
	\end{tabular}
	\par\end{center}

\subsection{SOSRI \label{subsec:SOSRI}}
\begin{center}
	\begin{tabular}{|c|c||c|c|}
		\hline
		Coefficient & Value & Coefficient & Value\tabularnewline
		\hline
		\hline
		$A_{2,1}^{(0)}$ & -0.04199224421316468 & $\alpha_{3}$ & 0.4736296532772559\tabularnewline
		\hline
		$A_{3,1}^{(0)}$ & 2.842612915017106 & $\alpha_{4}$ & 0.026404498125060714\tabularnewline
		\hline
		$A_{3,2}^{(0)}$ & -2.0527723684000727 & $c_{2}^{(0)}$ & -0.04199224421316468\tabularnewline
		\hline
		$A_{4,1}^{(0)}$ & 4.338237071435815 & $c_{3}^{(0)}$ & 0.7898405466170333\tabularnewline
		\hline
		$A_{4,2}^{(0)}$ & -2.8895936137439793 & $c_{4}^{(0)}$ & 3.7504010171562823\tabularnewline
		\hline
		$A_{4,3}^{(0)}$ & 2.3017575594644466 & $c_{1}^{(1)}$ & 0\tabularnewline
		\hline
		$A_{2,1}^{(1)}$ & 0.26204282091330466 & $c_{2}^{(1)}$ & 0.26204282091330466\tabularnewline
		\hline
		$A_{3,1}^{(1)}$ & 0.20903646383505375 & $c_{3}^{(1)}$ & 0.05879875232001766\tabularnewline
		\hline
		$A_{3,2}^{(1)}$ & -0.1502377115150361 & $c_{4}^{(1)}$ & 0.758661169101175\tabularnewline
		\hline
		$A_{4,1}^{(1)}$ & 0.05836595312746999 & $\beta_{1}^{(1)}$ & -1.8453464565104432\tabularnewline
		\hline
		$A_{4,2}^{(1)}$ & 0.6149440396332373 & $\beta_{2}^{(1)}$ & 2.688764531100726\tabularnewline
		\hline
		$A_{4,3}^{(1)}$ & 0.08535117634046772 & $\beta_{3}^{(1)}$ & -0.2523866501071323\tabularnewline
		\hline
		$B_{2,1}^{(0)}$ & -0.21641093549612528 & $\beta_{4}^{(1)}$ & 0.40896857551684956\tabularnewline
		\hline
		$B_{3,1}^{(0)}$ & 1.5336352863679572 & $\beta_{1}^{(2)}$ & 0.4969658141589478\tabularnewline
		\hline
		$B_{3,2}^{(0)}$ & 0.26066223492647056 & $\beta_{2}^{(2)}$ & -0.5771202869753592\tabularnewline
		\hline
		$B_{4,1}^{(0)}$ & -1.0536037558179159 & $\beta_{3}^{(2)}$ & -0.12919702470322217\tabularnewline
		\hline
		$B_{4,2}^{(0)}$ & 1.7015284721089472 & $\beta_{4}^{(2)}$ & 0.2093514975196336\tabularnewline
		\hline
		$B_{4,3}^{(0)}$ & -0.20725685784180017 & $\beta_{1}^{(3)}$ & 2.8453464565104425\tabularnewline
		\hline
		$B_{2,1}^{(1)}$ & -0.5119011827621657 & $\beta_{2}^{(3)}$ & -2.688764531100725\tabularnewline
		\hline
		$B_{3,1}^{(1)}$ & 2.67767339866713 & $\beta_{3}^{(3)}$ & 0.2523866501071322\tabularnewline
		\hline
		$B_{3,2}^{(1)}$ & -4.9395031322250995 & $\beta_{4}^{(3)}$ & -0.40896857551684945\tabularnewline
		\hline
		$B_{4,1}^{(1)}$ & 0.15580956238299215 & $\beta_{1}^{(4)}$ & 0.11522663875443433\tabularnewline
		\hline
		$B_{4,2}^{(1)}$ & 3.2361551006624674 & $\beta_{2}^{(4)}$ & -0.57877086147738\tabularnewline
		\hline
		$B_{4,3}^{(1)}$ & -1.4223118283355949 & $\beta_{3}^{(4)}$ & 0.2857851028163886\tabularnewline
		\hline
		$\alpha_{1}$ & 1.140099274172029 & $\beta_{4}^{(4)}$ & 0.17775911990655704\tabularnewline
		\hline
		$\alpha_{2}$ & -0.6401334255743456 &  & \tabularnewline
		\hline
	\end{tabular}
	\par\end{center}

\subsection{SOSRI2 \label{subsec:SOSRI2}}
\begin{center}
	\begin{tabular}{|c|c||c|c|}
		\hline
		Coefficient & Value & Coefficient & Value\tabularnewline
		\hline
		\hline
		$A_{2,1}^{(0)}$ & 0.13804532298278663 & $\alpha_{3}$ & 0.686995463807979\tabularnewline
		\hline
		$A_{3,1}^{(0)}$ & 0.5818361298250374 & $\alpha_{4}$ & -0.2911544680711602\tabularnewline
		\hline
		$A_{3,2}^{(0)}$ & 0.4181638701749618 & $c_{2}^{(0)}$ & 0.13804532298278663\tabularnewline
		\hline
		$A_{4,1}^{(0)}$ & 0.4670018408674211 & $c_{3}^{(0)}$ & 1\tabularnewline
		\hline
		$A_{4,2}^{(0)}$ & 0.8046204792187386 & $c_{4}^{(0)}$ & 1\tabularnewline
		\hline
		$A_{4,3}^{(0)}$ & -0.27162232008616016 & $c_{1}^{(1)}$ & 0\tabularnewline
		\hline
		$A_{2,1}^{(1)}$ & 0.45605532163856893 & $c_{2}^{(1)}$ & 0.45605532163856893\tabularnewline
		\hline
		$A_{3,1}^{(1)}$ & 0.7555807846451692 & $c_{3}^{(1)}$ & 1\tabularnewline
		\hline
		$A_{3,2}^{(1)}$ & 0.24441921535482677 & $c_{4}^{(1)}$ & 1\tabularnewline
		\hline
		$A_{4,1}^{(1)}$ & 0.6981181143266059 & $\beta_{1}^{(1)}$ & -0.45315689727309133\tabularnewline
		\hline
		$A_{4,2}^{(1)}$ & 0.3453277086024727 & $\beta_{2}^{(1)}$ & 0.8330937231303951\tabularnewline
		\hline
		$A_{4,3}^{(1)}$ & -0.04344582292908241 & $\beta_{3}^{(1)}$ & 0.3792843195533544\tabularnewline
		\hline
		$B_{2,1}^{(0)}$ & 0.08852381537667678 & $\beta_{4}^{(1)}$ & 0.24077885458934192\tabularnewline
		\hline
		$B_{3,1}^{(0)}$ & 1.0317752458971061 & $\beta_{1}^{(2)}$ & -0.4994383733810986\tabularnewline
		\hline
		$B_{3,2}^{(0)}$ & 0.4563552922077882 & $\beta_{2}^{(2)}$ & 0.9181786186154077\tabularnewline
		\hline
		$B_{4,1}^{(0)}$ & 1.73078280444124 & $\beta_{3}^{(2)}$ & -0.25613778661003145\tabularnewline
		\hline
		$B_{4,2}^{(0)}$ & -0.46089678470929774 & $\beta_{4}^{(2)}$ & -0.16260245862427797\tabularnewline
		\hline
		$B_{4,3}^{(0)}$ & -0.9637509618944188 & $\beta_{1}^{(3)}$ & 1.4531568972730915\tabularnewline
		\hline
		$B_{2,1}^{(1)}$ & 0.6753186815412179 & $\beta_{2}^{(3)}$ & -0.8330937231303933\tabularnewline
		\hline
		$B_{3,1}^{(1)}$ & -0.07452812525785148 & $\beta_{3}^{(3)}$ & -0.3792843195533583\tabularnewline
		\hline
		$B_{3,2}^{(1)}$ & -0.49783736486149366 & $\beta_{4}^{(3)}$ & -0.24077885458934023\tabularnewline
		\hline
		$B_{4,1}^{(1)}$ & -0.5591906709928903 & $\beta_{1}^{(4)}$ & -0.4976090683622265\tabularnewline
		\hline
		$B_{4,2}^{(1)}$ & 0.022696571806569924 & $\beta_{2}^{(4)}$ & 0.9148155835648892\tabularnewline
		\hline
		$B_{4,3}^{(1)}$ & -0.8984927888368557 & $\beta_{3}^{(4)}$ & -1.4102107084476505\tabularnewline
		\hline
		$\alpha_{1}$ & -0.15036858140642623 & $\beta_{4}^{(4)}$ & 0.9930041932449877\tabularnewline
		\hline
		$\alpha_{2}$ & 0.7545275856696072 &  & \tabularnewline
		\hline
	\end{tabular}
	\par\end{center}

\section{Appendix II: SRK Order Conditions }

\subsection{Order Conditions for Rößler-SRI Methods\label{subsec:Order-Conditions-for-SRI}}

The coefficients \\
$\left(A_{0},B_{0},\beta^{(i)},\alpha\right)$ must satisfy the following
order conditions to achieve order .5:

\begin{multicols}{3}
	\begin{enumerate}
		\item $\alpha^{T}e=1$
		\item $\beta^{(1)^{T}}e=1$
		\item $\beta^{(2)^{T}}e=0$
		\item $\beta^{(3)^{T}}e=0$
		\item $\beta^{(4)^{T}}e=0$
	\end{enumerate}
\end{multicols} additionally, for order 1:

\begin{multicols}{2}
	\begin{enumerate}
		\item $\beta^{(1)^{T}}B^{(1)}e=0$
		\item $\beta^{(2)^{T}}B^{(1)}e=1$
		\item $\beta^{(3)^{T}}B^{(1)}e=0$
		\item $\beta^{(4)^{T}}B^{(1)}e=0$
	\end{enumerate}
\end{multicols} and lastly for order 1.5:

\begin{multicols}{2}
	\begin{enumerate}
		\item $\alpha^{T}A^{(0)}e=\frac{1}{2}$
		\item $\alpha^{T}B^{(0)}e=1$
		\item $\alpha^{T}\left(B^{(0)}e\right)^{2}=\frac{3}{2}$
		\item $\beta^{(1)^{T}}A^{(1)}e=1$
		\item $\beta^{(2)^{T}}A^{(1)}e=0$
		\item $\beta^{(3)^{T}}A^{(1)}e=-1$
		\item $\beta^{(4)^{T}}A^{(1)}e=0$
		\item $\beta^{(1)^{T}}\left(B^{(1)}e\right)^{2}=1$
		\item $\beta^{(2)^{T}}\left(B^{(1)}e\right)^{2}=0$
		\item $\beta^{(3)^{T}}\left(B^{(1)}e\right)^{2}=-1$
		\item $\beta^{(4)^{T}}\left(B^{(1)}e\right)^{2}=2$
		\item $\beta^{(1)^{T}}\left(B^{(1)}\left(B^{(1)}e\right)\right)=0$
		\item $\beta^{(2)^{T}}\left(B^{(1)}\left(B^{(1)}e\right)\right)=0$
		\item $\beta^{(3)^{T}}\left(B^{(1)}\left(B^{(1)}e\right)\right)=0$
		\item $\beta^{(4)^{T}}\left(B^{(1)}\left(B^{(1)}e\right)\right)=1$
	\end{enumerate}
\end{multicols}
\begin{enumerate}
	\item [ ]
	\[
	16.\thinspace\thinspace\thinspace\frac{1}{2}\beta^{(1)^{T}}\left(A^{(1)}\left(B^{(0)}e\right)\right)+\frac{1}{3}\beta^{(3)^{T}}\left(A^{(1)}\left(B^{(0)}e\right)\right)=0
	\]
\end{enumerate}
where $f,g\in C^{1,2}(\mathcal{I}\times\mathbb{R}^{d},\mathbb{R}^{d})$,
$c^{(i)}=A^{(i)}e$, $e=(1,1,1,1)^{T}$ \cite{RN2707}.

\subsection{Order Conditions for Rößler-SRA Methods\label{subsec:Order-Conditions-for-SRA}}

The coefficients\\
$\left(A_{0},B_{0},\beta^{(i)},\alpha\right)$ must satisfy the conditions
for order 1:

\begin{multicols}{3}
	\begin{enumerate}
		\item $\alpha^{T}e=1$
		\item $\beta^{(1)^{T}}e=1$
		\item $\beta^{(2)^{T}}e=0$
	\end{enumerate}
\end{multicols}

and the additional conditions for order 1.5:

\begin{multicols}{3}
	\begin{enumerate}
		\item $\alpha^{T}B^{(0)}e=1$
		\item $\alpha^{T}A^{(0)}e=\frac{1}{2}$
		\item $\alpha^{T}\left(B^{(0)}e\right)^{2}=\frac{3}{2}$
		\item $\beta^{(1)^{T}}c^{(1)}=1$
		\item $\beta^{(2)^{T}}c^{(1)}=-1$
	\end{enumerate}
\end{multicols}

where $c^{(0)}=A^{(0)}e$ with $f\in C^{1,3}(\mathcal{I}\times\mathbb{R}^{d},\mathbb{R}^{d})$
and $g\in C^{1}(\mathcal{I},\mathbb{R}^{d})$ \cite{RN2707}.

\section{Appendix III: Derivation Details}

{\tiny{}
	\begin{eqnarray*}
		\left(I-\mu\Delta tA^{(0)}\right)H^{(0)} & = & U_{n}+\sigma\frac{I_{(1,0)}}{\Delta t}B^{(0)}\left(I-\sigma\sqrt{\Delta t}B^{(1)}\right)^{-1}\left(U_{n}+\mu\Delta tA^{(1)}H^{(0)}\right),\\
		\left(I-\mu\Delta tA^{(0)}\right)H^{(0)}-\left[\sigma\frac{I_{(1,0)}}{\Delta t}B^{(0)}\left(I-\sigma\sqrt{\Delta t}B^{(1)}\right)^{-1}\right]\mu\Delta tA^{(1)}H^{(0)} & = & U_{n}+\sigma\frac{I_{(1,0)}}{\Delta t}B^{(0)}\left(I-\sigma\sqrt{\Delta t}B^{(1)}\right)^{-1}U_{n}\\
		\left(I-\mu\Delta tA^{(0)}-\mu\Delta tA^{(1)}\sigma\frac{I_{(1,0)}}{\Delta t}B^{(0)}\left(I-\sigma\sqrt{\Delta t}B^{(1)}\right)^{-1}\right)H^{(0)} & = & \left(I+\sigma\frac{I_{(1,0)}}{\Delta t}B^{(0)}\left(I-\sigma\sqrt{\Delta t}B^{(1)}\right)^{-1}\right)U_{n}\\
		H^{(0)} & = & \left(I-\mu\Delta tA^{(0)}-\mu\sigma I_{(1,0)}A^{(1)}B^{(0)}\left(I-\sigma\sqrt{\Delta t}B^{(1)}\right)^{-1}\right)^{-1}\\
		&  & \left(I+\sigma\frac{I_{(1,0)}}{\Delta t}B^{(0)}\left(I-\sigma\sqrt{\Delta t}B^{(1)}\right)^{-1}\right)U_{n}
	\end{eqnarray*}
	\begin{eqnarray*}
		\left(I-\sigma\sqrt{\Delta t}B^{(1)}\right)H^{(1)} & = & U_{n}+\mu\Delta tA^{(1)}\left(I-\mu\Delta tA^{(0)}\right)^{-1}\left(U_{n}+\sigma\frac{I_{(1,0)}}{\Delta t}B^{(0)}H^{(1)}\right)\\
		\left(I-\sigma\sqrt{\Delta t}B^{(1)}-\mu\Delta tA^{(1)}\left(I-\mu\Delta tA^{(0)}\right)^{-1}\sigma\frac{I_{(1,0)}}{\Delta t}B^{(0)}\right)H^{(1)} & = & U_{n}+\mu\Delta tA^{(1)}\left(I-\mu\Delta tA^{(0)}\right)^{-1}U_{n}\\
		\left(I-\sigma\sqrt{\Delta t}B^{(1)}-\mu\Delta tA^{(1)}\left(I-\mu\Delta tA^{(0)}\right)^{-1}\sigma\frac{I_{(1,0)}}{\Delta t}B^{(0)}\right)H^{(1)} & = & \left(I+\mu\Delta tA^{(1)}\left(I-\mu\Delta tA^{(0)}\right)^{-1}\right)U_{n}\\
		H^{(1)} & = & \left(I-\sigma\sqrt{\Delta t}B^{(1)}-\mu\Delta tA^{(1)}\left(I-\mu\Delta tA^{(0)}\right)^{-1}\sigma\frac{I_{(1,0)}}{\Delta t}B^{(0)}\right)^{-1}\\
		&  & \left(I+\mu\Delta tA^{(1)}\left(I-\mu\Delta tA^{(0)}\right)^{-1}\right)U_{n}
	\end{eqnarray*}
}{\tiny \par}

{\tiny{}
	\begin{eqnarray*}
		U_{n+1} & = & U_{n}+\mu\Delta t\left(\alpha\cdot\left[\left(I-\mu\Delta tA^{(0)}-\mu\sigma I_{(1,0)}A^{(1)}B^{(0)}\left(I-\sigma\sqrt{\Delta t}B^{(1)}\right)^{-1}\right)^{-1}\left(I+\sigma\frac{I_{(1,0)}}{\Delta t}B^{(0)}\left(I-\sigma\sqrt{\Delta t}B^{(1)}\right)^{-1}\right)\right]U_{n}\right)\\
		&  & +\sigma I_{(1)}\left(\beta^{(1)}\cdot\left[\left(I-\sigma\sqrt{\Delta t}B^{(1)}-\mu\Delta tA^{(1)}\left(I-\mu\Delta tA^{(0)}\right)^{-1}\sigma\frac{I_{(1,0)}}{\Delta t}B^{(0)}\right)^{-1}\left(I+\mu\Delta tA^{(1)}\left(I-\mu\Delta tA^{(0)}\right)^{-1}\right)\right]U_{n}\right)\\
		&  & +\sigma\frac{I_{(1,1)}}{\sqrt{\Delta t}}\left(\beta^{(2)}\cdot\left[\left(I-\sigma\sqrt{\Delta t}B^{(1)}-\mu\Delta tA^{(1)}\left(I-\mu\Delta tA^{(0)}\right)^{-1}\sigma\frac{I_{(1,0)}}{\Delta t}B^{(0)}\right)^{-1}\left(I+\mu\Delta tA^{(1)}\left(I-\mu\Delta tA^{(0)}\right)^{-1}\right)\right]U_{n}\right)\\
		&  & +\sigma\frac{I_{(1,0)}}{\Delta t}\left(\beta^{(3)}\cdot\left[\left(I-\sigma\sqrt{\Delta t}B^{(1)}-\mu\Delta tA^{(1)}\left(I-\mu\Delta tA^{(0)}\right)^{-1}\sigma\frac{I_{(1,0)}}{\Delta t}B^{(0)}\right)^{-1}\left(I+\mu\Delta tA^{(1)}\left(I-\mu\Delta tA^{(0)}\right)^{-1}\right)\right]U_{n}\right)\\
		&  & +\sigma\frac{I_{(1,1,1)}}{\Delta t}\left(\beta^{(4)}\cdot\left[\left(I-\sigma\sqrt{\Delta t}B^{(1)}-\mu\Delta tA^{(1)}\left(I-\mu\Delta tA^{(0)}\right)^{-1}\sigma\frac{I_{(1,0)}}{\Delta t}B^{(0)}\right)^{-1}\left(I+\mu\Delta tA^{(1)}\left(I-\mu\Delta tA^{(0)}\right)^{-1}\right)\right]U_{n}\right)
	\end{eqnarray*}
} Thus we substitute in the Wiktorsson approximations
\begin{align*}
I_{(i,i)} & =\frac{1}{2}\left(\Delta W^{2}-h\right)\\
I_{(i,i,i)} & =\frac{1}{6}\left(\Delta W^{3}-3h\Delta W\right)\\
I_{(i.0)} & =\frac{1}{2}h\left(\Delta W+\frac{1}{\sqrt{3}}\Delta Z\right)
\end{align*}
where $\Delta Z\sim N(0,h)$ is independent of $\Delta W\sim N(0,h)$.
By the properties of the normal distribution, we have that

\[
E\left[\left(\Delta W\right)^{n}\right]=0
\]
for any odd $n$ and
\begin{eqnarray*}
	E\left[\left(\Delta W\right)^{2}\right] & = & h\\
	E\left[\left(\Delta W\right)^{4}\right] & = & 3h^{2}\\
	E\left[\left(\Delta W\right)^{6}\right] & = & 15h^{3}\\
	E\left[\left(\Delta W\right)^{8}\right] & = & 105h^{4},
\end{eqnarray*}
and similarly for $\Delta Z$.

\section*{Acknowledgments}

We would like to thank the members of JuliaDiffEq, specifically David Widmann (@devmotion),
Yingbo Ma (@YingboMa) and (@dextorious) for their contributions to the ecosystem. Their efforts
have helped make the development of efficient implementations possible.

\bibliographystyle{siamplain}
\bibliography{references}

\begin{thebibliography}{10}

\bibitem{RN3800}
{\sc A.~Abdulle and S.~Cirilli}, {\em S-rock: Chebyshev methods for stiff
  stochastic differential equations}, SIAM Journal on Scientific Computing, 30
  (2008), pp.~997--1014, \url{https://doi.org/10.1137/070679375},
  \url{https://doi.org/10.1137/070679375}.

\bibitem{RN3515}
{\sc D.~I.~K. Ahmadia and A.~J.}, {\em Optimal stability polynomials for
  numerical integration of initial value problems}, CAMCOS, 7 (2012),
  pp.~247--271, \url{https://doi.org/10.2140/camcos.2012.7.247}.

\bibitem{RN3798}
{\sc E.~G. Birgin and J.~M. Martinez}, {\em Improving ultimate convergence of
  an augmented lagrangian method}, Optimization Methods and Software, 23
  (2008), pp.~177--195, \url{https://doi.org/10.1080/10556780701577730},
  \url{https://doi.org/10.1080/10556780701577730}.

\bibitem{RN3791}
{\sc K.~Burrage and J.~C. Butcher}, {\em Non-linear stability of a general
  class of differential equation methods}, BIT Numerical Mathematics, 20
  (1980), pp.~185--203, \url{https://doi.org/10.1007/BF01933191},
  \url{https://doi.org/10.1007/BF01933191}.

\bibitem{RN3516}
{\sc J.~C. Butcher}, {\em A history of runge-kutta methods}, Applied Numerical
  Mathematics, 20 (1996), pp.~247--260,
  \url{https://doi.org/http://dx.doi.org/10.1016/0168-9274(95)00108-5},
  \url{http://www.sciencedirect.com/science/article/pii/0168927495001085}.

\bibitem{RN3519}
{\sc J.~C. Butcher}, {\em Numerical methods for ordinary differential equations
  in the 20th century}, Journal of Computational and Applied Mathematics, 125
  (2000), pp.~1--29,
  \url{https://doi.org/http://doi.org/10.1016/S0377-0427(00)00455-6},
  \url{http://www.sciencedirect.com/science/article/pii/S0377042700004556}.

\bibitem{RN3797}
{\sc A.~Conn, N.~Gould, and P.~Toint}, {\em A globally convergent augmented
  lagrangian algorithm for optimization with general constraints and simple
  bounds}, SIAM Journal on Numerical Analysis, 28 (1991), pp.~545--572,
  \url{https://doi.org/10.1137/0728030}, \url{https://doi.org/10.1137/0728030}.

\bibitem{RN3513}
{\sc J.~R. Dormand and P.~J. Prince}, {\em A family of embedded runge-kutta
  formulae}, Journal of Computational and Applied Mathematics, 6 (1980),
  pp.~19--26,
  \url{https://doi.org/http://dx.doi.org/10.1016/0771-050X(80)90013-3},
  \url{http://www.sciencedirect.com/science/article/pii/0771050X80900133}.

\bibitem{RN3795}
{\sc I.~Dunning, J.~Huchette, and M.~Lubin}, {\em Jump: A modeling language for
  mathematical optimization}, SIAM Review, 59 (2017), pp.~295--320,
  \url{https://doi.org/10.1137/15M1020575},
  \url{https://doi.org/10.1137/15M1020575}.

\bibitem{RN3802}
{\sc H.~Gilsing and T.~Shardlow}, {\em Sdelab: A package for solving stochastic
  differential equations in matlab}, Journal of Computational and Applied
  Mathematics, 205 (2007), pp.~1002--1018,
  \url{https://doi.org/https://doi.org/10.1016/j.cam.2006.05.037},
  \url{http://www.sciencedirect.com/science/article/pii/S0377042706004195}.

\bibitem{RN3353}
{\sc E.~Hairer, S.~P. Nørsett, and G.~Wanner}, {\em Solving ordinary
  differential equations I : nonstiff problems}, Springer series in
  computational mathematics,, Springer, Heidelberg ; London, 2nd rev.~ed.,
  2009.

\bibitem{RN3790}
{\sc E.~Hairer and G.~Wanner}, {\em Solving Ordinary Differential Equations II
  - Stiff and Differential-Algebraic Problems}, Springer, 1991.

\bibitem{RN3794}
{\sc E.~Hairer and G.~Wanner}, {\em Stiff differential equations solved by
  radau methods}, Journal of Computational and Applied Mathematics, 111 (1999),
  pp.~93--111,
  \url{https://doi.org/https://doi.org/10.1016/S0377-0427(99)00134-X},
  \url{http://www.sciencedirect.com/science/article/pii/S037704279900134X}.

\bibitem{RN3358}
{\sc T.~Hong, K.~Watanabe, C.~H. Ta, A.~Villarreal-Ponce, Q.~Nie, and X.~Dai},
  {\em An ovol2-zeb1 mutual inhibitory circuit governs bidirectional and
  multi-step transition between epithelial and mesenchymal states}, PLoS Comput
  Biol, 11 (2015), p.~e1004569,
  \url{https://doi.org/10.1371/journal.pcbi.1004569},
  \url{http://dx.doi.org/10.1371%2Fjournal.pcbi.1004569}.

\bibitem{RN3793}
{\sc M.~E. Hosea and L.~F. Shampine}, {\em Analysis and implementation of
  tr-bdf2}, Applied Numerical Mathematics, 20 (1996), pp.~21--37,
  \url{https://doi.org/https://doi.org/10.1016/0168-9274(95)00115-8},
  \url{http://www.sciencedirect.com/science/article/pii/0168927495001158}.

\bibitem{RN3804}
{\sc A.~Janicki, A.~Izydorczyk, and P.~Gradalski}, {\em Computer Simulation of
  Stochastic Models with SDE-Solver Software Package}, Springer Berlin
  Heidelberg, Berlin, Heidelberg, 2003, pp.~361--370,
  \url{https://doi.org/10.1007/3-540-44860-8_37},
  \url{https://doi.org/10.1007/3-540-44860-8_37}.

\bibitem{RN3788}
{\sc S.~G. Johnson}, {\em The nlopt nonlinear-optimization package},
  \url{http://ab-initio.mit.edu/nlopt}.

\bibitem{RN3792}
{\sc C.~A. Kennedy and M.~H. Carpenter}, {\em Additive runge-kutta schemes for
  convection-diffusion-reaction equations}, Applied Numerical Mathematics, 44
  (2003), pp.~139--181,
  \url{https://doi.org/https://doi.org/10.1016/S0168-9274(02)00138-1},
  \url{http://www.sciencedirect.com/science/article/pii/S0168927402001381}.

\bibitem{RN3815}
{\sc P.~Kloeden and A.~Neuenkirch}, {\em Convergence of numerical methods for
  stochastic differential equations in mathematical finance}, 2012,
  \url{https://doi.org/10.1142/9789814436434_0002}.

\bibitem{RN3169}
{\sc P.~E. Kloeden and E.~Platen}, {\em Numerical Solution of Stochastic
  Differential Equations}, Springer Berlin Heidelberg, 2011,
  \url{https://books.google.com/books?id=BCvtssom1CMC}.

\bibitem{RN3801}
{\sc Y.~Komori and K.~Burrage}, {\em Weak second order s-rock methods for
  stratonovich stochastic differential equations}, Journal of Computational and
  Applied Mathematics, 236 (2012), pp.~2895--2908,
  \url{https://doi.org/https://doi.org/10.1016/j.cam.2012.01.033},
  \url{http://www.sciencedirect.com/science/article/pii/S0377042712000441}.

\bibitem{RN3799}
{\sc Y.~Komori and K.~Burrage}, {\em Strong first order s-rock methods for
  stochastic differential equations}, Journal of Computational and Applied
  Mathematics, 242 (2013), pp.~261--274,
  \url{https://doi.org/https://doi.org/10.1016/j.cam.2012.10.026},
  \url{http://www.sciencedirect.com/science/article/pii/S0377042712004669}.

\bibitem{RN3512}
{\sc F.~S. Lawrence}, {\em Some practical runge-kutta formulas}, Math. Comput.,
  46 (1986), pp.~135--150, \url{https://doi.org/10.2307/2008219}.

\bibitem{RN3514}
{\sc J.~Lawson}, {\em An order five runge-kutta process with extended region of
  stability}, SIAM Journal on Numerical Analysis, 3 (1966), pp.~593--597,
  \url{https://doi.org/10.1137/0703051},
  \url{http://dx.doi.org/10.1137/0703051}.

\bibitem{RN3811}
{\sc T.~Li, A.~Abdulle, and W.~E}, {\em Effectiveness of implicit methods for
  stiff stochastic differential equations}, Commun. Comput. Phys, 3 (2008),
  pp.~295--307.

\bibitem{RN3806}
{\sc X.~Mao}, {\em The truncated euler–maruyama method for stochastic
  differential equations}, Journal of Computational and Applied Mathematics,
  290 (2015), pp.~370--384,
  \url{https://doi.org/https://doi.org/10.1016/j.cam.2015.06.002},
  \url{http://www.sciencedirect.com/science/article/pii/S0377042715003210}.

\bibitem{RN3525}
{\sc J.~K. Moller and H.~Madsen}, {\em From state dependent diffusion to
  constant diffusion in stochastic differential equations by the lamperti
  transform}, report, Technical University of Denmark, DTU Informatics,
  Building 321, 2010.

\bibitem{RN3789}
{\sc M.~J.~D. Powell}, {\em A direct search optimization method that models the
  objective and constraint functions by linear interpolation}, in Advances in
  Optimization and Numerical Analysis, Proceedings of the 6th Workshop on
  Optimization and Numerical Analysis, Oaxaca, Mexico, S.~Gomez and J.-P.
  Hennart, eds., vol.~275, Kluwer Academic Publishers, pp.~51--67,
  \url{https://doi.org/citeulike-article-id:6904064},
  \url{http://www.ams.org/mathscinet-getitem?mr=95d:90075}.

\bibitem{RN3787}
{\sc C.~Rackauckas and Q.~Nie}, {\em Adaptive methods for stochastic
  differential equations via natural embeddings and rejection sampling with
  memory}, Discrete and Continuous Dynamical Systems - Series B, 22 (2016),
  pp.~2731--2761, \url{https://doi.org/10.3934/dcdsb.2017133},
  \url{http://aimsciences.org//article/id/5354a27a-e5be-4c40-9d7a-e918853b56b7}.

\bibitem{RN3784}
{\sc C.~Rackauckas and Q.~Nie}, {\em Differentialequations.jl - a performant
  and feature-rich ecosystem for solving differential equations in julia},
  Journal of Open Research Software, 5 (2017), p.~15,
  \url{https://doi.org/http://doi.org/10.5334/jors.151}.

\bibitem{RN3805}
{\sc C.~Rackauckas and Q.~Nie}, {\em Mean-independent noise control of cell
  fates via intermediate states}, iScience, Accepted (2018).

\bibitem{RN2707}
{\sc A.~Rossler}, {\em Runge kutta methods for the strong approximation of
  solutions of stochastic differential equations}, SIAM Journal on Numerical
  Analysis, 48 (2010), pp.~922--952, \url{https://doi.org/10.1137/09076636x}.

\bibitem{RN3796}
{\sc T.~P. Runarsson and Y.~Xin}, {\em Search biases in constrained
  evolutionary optimization}, IEEE Transactions on Systems, Man, and
  Cybernetics, Part C (Applications and Reviews), 35 (2005), pp.~233--243,
  \url{https://doi.org/10.1109/TSMCC.2004.841906}.

\bibitem{RN3803}
{\sc T.~Schaffter}, {\em From genes to organisms: Bioinformatics System Models
  and Software}, thesis, 2014.

\bibitem{RN3527}
{\sc L.~F. Shampine}, {\em Stiffness and nonstiff differential equation
  solvers, ii: Detecting stiffness with runge-kutta methods}, ACM Trans. Math.
  Softw., 3 (1977), pp.~44--53, \url{https://doi.org/10.1145/355719.355722}.

\bibitem{RN3526}
{\sc L.~F. Shampine and K.~L. Hiebert}, {\em Detecting stiffness with the
  fehlberg (4, 5) formulas}, Computers \& Mathematics with Applications, 3
  (1977), pp.~41--46,
  \url{https://doi.org/http://dx.doi.org/10.1016/0898-1221(77)90112-2},
  \url{http://www.sciencedirect.com/science/article/pii/0898122177901122}.

\bibitem{RN3518}
{\sc W.~H.~E. Sharp, D.~J. Higham, B.~Owren, and P.~W.}, {\em A survey of the
  explicit runge-kutta method},  (1995).

\bibitem{RN3807}
{\sc G.~Soderlind and L.~Wang}, {\em Evaluating numerical ode/dae methods,
  algorithms and software}, Journal of Computational and Applied Mathematics,
  185 (2006), pp.~244--260,
  \url{https://doi.org/https://doi.org/10.1016/j.cam.2005.03.009},
  \url{http://www.sciencedirect.com/science/article/pii/S0377042705001135}.

\bibitem{RN3523}
{\sc C.~Tsitouras}, {\em Runge-kutta pairs of order 5(4) satisfying only the
  first column simplifying assumption}, Computers \& Mathematics with
  Applications, 62 (2011), pp.~770--775,
  \url{https://doi.org/http://doi.org/10.1016/j.camwa.2011.06.002},
  \url{http://www.sciencedirect.com/science/article/pii/S0898122111004706}.

\bibitem{RN3517}
{\sc P.~J. van~der Houwen}, {\em Explicit runge-kutta formulas with increased
  stability boundaries}, Numerische Mathematik, 20 (1972), pp.~149--164,
  \url{https://doi.org/10.1007/BF01404404},
  \url{http://dx.doi.org/10.1007/BF01404404}.

\bibitem{RN3814}
{\sc B.~Wang and Q.~Qi}, {\em Modeling the lake eutrophication stochastic
  ecosystem and the research of its stability}, Mathematical Biosciences,
  \url{https://doi.org/https://doi.org/10.1016/j.mbs.2018.03.019},
  \url{https://www.sciencedirect.com/science/article/pii/S0025556418301780}.

\bibitem{RN3816}
{\sc P.~Wang and Y.~Li}, {\em Split-step forward methods for stochastic
  differential equations}, Journal of Computational and Applied Mathematics,
  233 (2010), pp.~2641--2651,
  \url{https://doi.org/https://doi.org/10.1016/j.cam.2009.11.010},
  \url{http://www.sciencedirect.com/science/article/pii/S0377042709007419}.

\bibitem{RN3175}
{\sc M.~Wiktorsson}, {\em Joint characteristic function and simultaneous
  simulation of iterated ito integrals for multiple independent brownian
  motions}, The Annals of Applied Probability,  (2001), pp.~470--487,
  \url{https://doi.org/10.1214/aoap/1015345301},
  \url{http://projecteuclid.org/euclid.aoap/1015345301}.

\end{thebibliography}
\end{document}